\documentclass{article}
\usepackage{amssymb}
\usepackage{amsmath,color}
\def\thebibliography#1{\section*{References}\list
  {[\arabic{enumi}]}{\settowidth\labelwidth{[#1]}\leftmargin\labelwidth
    \advance\leftmargin\labelsep
    \usecounter{enumi}}
    \def\newblock{\hskip .11em plus .33em minus -.07em}
    \sloppy
    \sfcode`\.=1000\relax}

\newcommand{\refbook}[3]{{\sc #1}{\em\ #2}{\ #3}}
\newcommand{\refer}[5]{{\sc #1}{\ #2}{\em\ #3}{\bf\ #4}{\ #5}}

\newtheorem{lem}{Lemma}[section]
\newtheorem{cor}[lem]{Corollary}
\newtheorem{teo}[lem]{Theorem}
\newtheorem{os}[lem]{Remark}

\newtheorem{prop}[lem]{Proposition}

\newcommand{\qed}{\thinspace\null\nobreak\hfill\hbox{\vbox{\kern-.2pt\hrule
 height.2pt depth.2pt\kern-.2pt\kern-.2pt \hbox to2.5mm{\kern-.2pt\vrule
 width.4pt \kern-.2pt\raise2.5mm\vbox to.2pt{}\lower0pt\vtop
 to.2pt{}\hfil\kern-.2pt \vrule
 width.4pt \kern-.2pt}\kern-.2pt\kern-.2pt\hrule height.2pt depth.2pt
 \kern-.2pt}}\par\medbreak}

\newcommand{\R}{\mathbb{R}}

\newcommand{\C}{\mathbb{C}}

\newcommand{\Rp}{\textrm{\emph{Re}\,}}

\newcommand{\eps}{\varepsilon}

\newcommand{\ov}{\overline}

\newcommand{\ds}{\displaystyle}

\newcommand{\aver}[1]{-\hskip-0.38cm\int_{#1}}

\date{}
\begin{document}

\title{ $L^p$ estimates for the Caffarelli-Silvestre extension operators}
\author{G. Metafune \thanks{Dipartimento di Matematica e Fisica ``Ennio De Giorgi'', Universit\`a del Salento, C.P.193, 73100, Lecce, Italy.
-mail:  giorgio.metafune@unisalento.it}\qquad L. Negro \thanks{Dipartimento di Matematica e Fisica  ``Ennio De
Giorgi'', Universit\`a del Salento, C.P.193, 73100, Lecce, Italy. email: luigi.negro@unisalento.it} \qquad C. Spina \thanks{Dipartimento di Matematica e Fisica``Ennio De Giorgi'', Universit\`a del Salento, C.P.193, 73100, Lecce, Italy.
e-mail:  chiara.spina@unisalento.it}}

\maketitle
\begin{abstract}
\noindent 
We study elliptic and parabolic problems governed by the  singular elliptic   operators 
 \begin{equation*}
\mathcal L =\Delta_{x} +D_{yy}+\frac{c}{y}D_y  -\frac{b}{y^2}
\end{equation*}
in the half-space $\R^{N+1}_+=\{(x,y): x \in \R^N, y>0\}$.
 
\bigskip\noindent
Mathematics subject classification (2010): 47D07, 35J70.
\par

\noindent Keywords: elliptic operators, discontinuous coefficients, kernel estimates, maximal regularity.
\end{abstract}

\section{Introduction}
 In this paper we study solvability and regularity of elliptic and parabolic problems associated to the  degenerate   operators 
\begin{equation*} \label{defL}
\mathcal L =\Delta_{x} +D_{yy}+\frac{c}{y}D_y  -\frac{b}{y^2} \quad {\rm and}\quad D_t- \mathcal L
\end{equation*}
in the half-space $\R^{N+1}_+=\{(x,y): x \in \R^N, y>0\}$ or  $(0, \infty) \times \R^{N+1}_+$.

Here $b,\ c$ are constant real coefficients and we use 
$
L_y=D_{yy}+\frac{c}{y}D_y  -\frac{b}{y^2}.
$
 Note that singularities in the lower order terms appear when either $b$ or $c$ is different from 0. 

The operators $\Delta_x$, $L_y$ commute  and the whole operator $\mathcal L$ satisfies the scaling property $I_s^{-1}\mathcal L I_s=s^2\mathcal L$, if $I_s u(x,y)=u(sx,sy)$.

When $b=0$, then $L_y$ is a Bessel operator (we shall denote it by $B_y$) and both $\mathcal L=\Delta_x+B_y$ and $D_t-\mathcal L$ play a major role in the investigation of the fractional powers $(-\Delta_x)^s$ and  $(D_t-\Delta_x)^s$, $s=(1-c)/2$, through the  ``extension procedure" of Caffarelli and Silvestre \cite{CS}, after the pioneering work by Muckenhoupt and Stein, \cite{MS}. For this reason, ${\mathcal L}$ and $D_t-\mathcal L$ are named the ``extension operators".
We refer the reader to \cite[Section 10]{Garofalo} for an exposition of  the theory in the language of semigroups and for the references to the wide literature on the extension problem, both in the elliptic and parabolic case.

Here we study unique solvability of the problems $\lambda u-\mathcal L u=f$ and $D_t v -\mathcal L v=g$ in $L^p$ spaces under appropriate  boundary conditions, and initial conditions in the parabolic case, together with the regularity of $u,v$. In the language of semigroup theory, we prove that $\mathcal L$ generates an analytic semigroup, characterize its domain and show that it has maximal regularity, which means that both $D_t v$ and $\mathcal L v$ have the same regularity as $g$.


Both the domains of $\Delta_x$ and $L_y$ are known, in their corresponding $L^p$-spaces. Clearly $D(\Delta_x)=W^{2,p}(\R^N)$, $1<p<\infty$. However $D(L_y)\subset L^p(0, \infty)$ is more delicate, the boundary conditions and the regularity up to $y=0$ depend on the coefficients $c,b$. We shall devote Sections 2, 3, 4  to a careful study of the 1d operator $L_y$, starting form the Bessel case where $b=0$, this last both under Dirichlet and Neumann boundary conditions at $y=0$. We shall provide  in both cases the description of the domain,  pointwise estimates for the heat kernel and  its gradient. The general case is reduced to the Bessel one, through a change of variables.
We study $L_y$ also in weighted spaces $L^p((0,\infty), y^mdy)$; the cases $m=0$ and $m=c$ are the most important: the first corresponds to the Lebesgue measure, the second to the symmetrizing one. However we need general $m$ also for technical reasons. This makes the exposition slightly heavier, but it is unavoidable in our approach. Not all the results in these sections are completely new. A description of the domain under Dirichlet boundary conditions is in \cite{met-soba-spi3} but here we have more precise results; we are not aware of a similar description in the case of Neumann boundary conditions for $B_y$. The heat kernel is known among probabilists but a purely analytic derivation can be found in \cite{Garofalo1} for Neumann boundary conditions. Here we prefer to give an analytic proof in both cases and provide  manageable and precise estimates.

The elliptic operator $\mathcal L$ is studied through estimates like 
\begin{equation} \label{closedness}
\|\Delta_x u\|_p+\|L_y u\|_p \le C\| \mathcal L u\|_p
\end{equation}
where the $L^p$ norms are taken over $\R_+^{N+1}$. This kind of estimates are quite natural in this context but not easy to prove. Of course they are equivalent to $\|D_{x_ix_j}u\|_p \le C\| \mathcal Lu\|_p$, by the Calder\'{o}n-Zygmund inequalities in the $x$-variables, and can be restated by saying that $\mathcal L$ is closed on $D(\Delta_x) \cap D(L_y)$ or that $\Delta_x  \mathcal L^{-1}$ is bounded. Note that the weaker inequality \eqref{closedness} with $\| \mathcal Lu\|_p+\|u\|_p$ on the right hand side implies the stronger one as stated, by the scaling properties of the  operators involved.

 Estimates for $D_{yy} u$ follow if and only if they hold for the one dimensional operator  $L_y$ but those for the mixed derivatives $D_{x_i y} u$ are more subtle. They are certainly true when  $D_{yy}L_y^{-1}$ is bounded, by Calder\'{o}n-Zygmund with respect to all $x,y$ variables, but we shall prove that they hold if (and only if) $D_y(I-L_y)^{-1}$ is bounded, which was quite unexpected for us.

Let us explain how to obtain \eqref{closedness} when $p=2$ and introduce our approach for general $p$. Assuming that $\Delta_x u+L_yu=f$ and taking the Fourier transform with respect to $x$ (with covariable $\xi$) we obtain $-|\xi|^2 \hat u(\xi,y)+L_y \hat u(\xi,y)=\hat f(\xi,y)$ and then $|\xi|^2 \hat u(\xi,y)=-|\xi|^2 (|\xi|^2-L_y)^{-1}\hat f (\xi,y)$. Assuming that $L_y$ generates a bounded semigroup in $L^2(0, \infty)$, then $|\xi|^2\|(|\xi|^2-L_y)^{-1}\| \le C$ and
$$
\int_0^\infty |\xi|^4 |\hat u (\xi,y)|^2dy \le C^2 \int_0^\infty |\hat f(\xi,y)|^2dy
$$
which gives, after integration with respect to $\xi$ and Plancherel equality,
$$
\|\Delta_x u\|_2 =\||\xi|^2 \hat u \|_2 \le C\|f\|_2.
$$
When $p \neq 2$ and denoting by ${\cal F}$ the Fourier transform with respect to $x$  we get, formally, 
$$
\Delta_x  \mathcal  L^{-1}=-{\cal F}^{-1}\left (|\xi|^2(|\xi|^2-L_y)^{-1} \right) {\cal F}
$$
and the boundedness of $\Delta_x  \mathcal L^{-1}$ is equivalent to say that the operator valued map $\xi \in \R^N \to M(\xi)=|\xi|^2(|\xi|^2-L_y)^{-1} \in B(L^p(0,\infty))$ is a bounded Fourier multiplier in $L^p(\R^N; L^p(0,\infty))=L^p(\R_+^{N+1})$.
Here we use a vector valued Mikhlin multiplier theorem which relies on the $\mathcal R$-boundedness of the family $M(\xi)$ and its derivatives, which we deduce from heat kernel estimates.
We use a similar strategy for $\nabla_x D_y  \mathcal L^{-1}$ which this time rests on estimates for the gradient of the heat kernel of $L_y$.

It is important to note that the closedness of $\Delta_x+L_y$ on the intersection of the corresponding domains does not follow from general results. In fact, $e^{tL_y}$ is not contractive and does not admit Gaussian estimates, except for special cases; moreover it is bounded in $L^p(0,\infty)$ only for certain intervals of $p$ depending on the coefficients $c,b$.

The strategy for proving  the parabolic estimates
$$
\|D_t v\|_p+\|\mathcal L v\|_p \le C\|(D_t- \mathcal L) v\|_p
$$
($L^p$ norms on $(0,\infty)\times \R^{N+1}_+$), is similar after taking the Fourier transform with respect to $t$.

Both the elliptic and parabolic estimates rely on a vector valued Mikhlin multiplier theorem and share the name ``maximal regularity" even though this term is often restricted to the parabolic case.

The functional analytic approach for  maximal regularity is widely described in \cite{KW} and in the new books \cite{WeisBook1}, \cite{WeisBook2}. The whole theory relies on a deep interplay between harmonic analysis and structure theory of Banach spaces but largely simplifies when the underlying Banach spaces are $L^p$ spaces, by using classical square function estimates. This last approach has been employed extensively in \cite{DHP}, showing that uniformly parabolic operators have maximal regularity, under very general boundary conditions. 

We deduce the boundedness of vector valued multipliers by the  $\mathcal R$-boundedness of  a family of integral operators, which we prove  through  an extrapolation result in \cite{auscher1} which involves a family of Muckenhoupt weighted estimates.  Here we adopt the same strategy as T. A.  Bui, see \cite{bui}, in the case of Schr\"odinger operators with inverse square potentials. Section 7  is really the core of the paper, while Section 6 contains all relevant definitions and results for the subsequent proofs.

We work in  $L^p(\R^{N+1}_+, y^m dx dy)$ not just for the sake of generality but because our proof relies on weighted estimates: we are unable to obtain the result just fixing the Lebesgue measure or the symmetrizing one $y^c dx dy$ but we have to work simultaneously in different homogeneous spaces.   

As an application of our results, in Section 9 we deduce Rellich inequalities for $\mathcal L=\Delta_x+L_y$ by the analogous results for the one dimensional operator $L_y$, using the closedness of $\mathcal L$ on the intersection of the domains of $\Delta_x$ and $L_y$.

\bigskip
\noindent\textbf{Notation.} For $N \ge 0$, $\R^{N+1}_+=\{(x,y): x \in \R^N, y>0\}$. For $m \in \R$ we consider the measure $d\mu_m =y^m dx dy $ in $\R^{N+1}_+$. We write $L^p_m(\R^{N+1})$ for  $L^p(\R_+^{N+1}; y^m dx dy)$ and often only $L^p_m$ when $\R^{N+1}_+$ is understood. Similarly $W^{k,p}_m(\R^{N+1}_+)=\{u \in L^p_m(\R^{N+1}_+): \partial^\alpha u \in  L^p_m(\R^{N+1}_+) \quad |\alpha| \le k\}$. We use often $W^{k,p}_m$ thus omitting $\R^{N+1}_+$ and $W^{k,p}_{0,m}$ for the closure of $C_c^\infty (\R^{N+1}_+)$ in $W^{k,p}_m$ and we use $H^k_m$ for $W^{k,2}_m$.

\bigskip

\noindent\textbf{Acknowledgements.} The authors thank S. Fornaro, N. Garofalo,  D. Pallara and V. Vespri  for several comments on a previous version of the manuscript.

\section{Bessel operators in $1d$}

In this section we state and prove the main properties of the degenerate operator
$$
B=D_{yy}+\frac{c}{y}D_y=y^{-c}D_y \left (y^c D_y\right )
$$ on the half line $\R_+=]0, \infty[$
needed for our purposes.

\subsection{Weighted $L^2$ spaces and Bessel operators}
We use the Sobolev spaces defined in Appendix B, for $p=2$ and $N=0$. According to the above notation, for  $c \in \R$ we use $L^2_c=\{u: \R_+ \to \C: \int_0^\infty |u(y)|^2 y^c dy <\infty\}$, $H^1_c=\{u \in L^2_c, u' \in L^2_c\}$, where $u'$ is understood as a distribution in the open interval $]0, \infty[$. Both $L^2_c$ and $H^1_c$ are Hilbert spaces under their canonical inner products; moreover $C_c^\infty (0, \infty)$ is contained in both and dense in $L^2_c$. We denote by $H^1_{0,c}$ the closure of $C_c^\infty (0, \infty)$ in $H^1_c$. We need the following properties proved in greater generality in Appendix B.

\begin{lem}\label{capacity}
\begin{itemize}
\item[(i)] If $|c |\geq 1$, then $H^1_{0,c}=H^1_c$. When $c \le -1$, then $\ds \lim_{y \to 0}u(y)=0$ for every $u \in H^1_c$.
\item[(ii)] If $|c|  <1$ and $u \in H^1_c$, then $\ds\lim_{y \to 0}u(y)=\ell \in \C$. Moreover, $\ell=0$ if and only if $u \in H^1_{0, c}$.
\end{itemize}
\end{lem}

 $B$ is associated to the  symmetric form  in $L^2_c$
 
\begin{align*}
\mathfrak{a}(u,v)
&:=
\int_0^\infty D_y u  D_y \overline{v}\,y^c dy=\int_0^\infty ( Bu)\, \overline{v}\, y^c dy.
\end{align*}
For any $c \in \R$ we may consider $H^1_{0,c}$ as domain of the form and, accordingly, define the Bessel  operator with Dirichlet boundary conditions $B^d$ by
\begin{equation} \label{BesselD}
D(B^d)=\{u \in H^1_{0,c}: \exists  f \in L^2_c \ {\rm such\ that}\  \mathfrak{a}(u,v)=\int_0^\infty f \overline{v}y^c\, dy\ {\rm for\ every}\ v\in H^1_{0,c}\}, \quad B^du=-f
\end{equation}
Similarly, by considering $H^1_c$ we obtain the Bessel operator with Neumann Boundary conditions $B^n$ defined as
\begin{equation} \label{BesselN}
D(B^n)=\{u \in H^1_{c}: \exists  f \in L^2_c \ {\rm such\ that}\  \mathfrak{a}(u,v)=\int_0^\infty f \overline{v}y^c\, dy\ {\rm for\ every}\ v\in H^1_{c}\}, \quad B^nu=-f
\end{equation}
$B^d, B^n$ are non-positive self-adjoint operators and $ u \in H^2_{loc}(0, \infty)$ with $B^du=B^nu=u_{yy}+\frac{c}{y}u_y$ for $y>0$,  if $u \in D(B^d)$ or $u \in D(B^n)$, by standard arguments.

\begin{lem} \label{Neumann} If $c>-1$ and $u \in D(B^n)$, then $\ds\lim_{y \to 0}y^cu'(y)=0$.
\end{lem}
{\sc Proof. } By  assumption for $v \in H^1_c$
\begin{align*}
\int_0^\infty u_y v_y y^cdy&=-\int_0^\infty (B^n u) v y^c dy=-\lim_{\eps \to 0}\int_\eps^\infty \frac{d}{dy} (y^c u_y)v dy\\
&=\lim_{\eps \to 0}\left (\int_\eps^\infty  u_y v_y y^c dy -\eps^c u_y(\eps)v(\eps)\right ) =\int_0^\infty u_y v_y y^c dy -\lim_{\eps \to 0} \eps^c u_y(\eps)v(\eps).
\end{align*}

Choosing $v \equiv 1$ near $0$, which is possible since $c>-1$, we get the result. \qed

\medskip
Observe that: 
\begin{itemize}
\item  when $|c| \ge 1$ then $B^d=B^n$  and, when $c \le -1$, $u(0)=0$ for every $u \in D(B^d)$ by Lemma \ref{capacity} (i);
\item when $|c|<1$ then  $B^d$ and $B^n$ are different and $u \in D(B^d)$ fulfils $u(0)=0$, by Lemma \ref{capacity} (ii).
\end{itemize}
%
%

Even though $B^d$ and $B^n$ are defined for every $c \in \R$, we shall use  $B^d$  when $c<1$ and $B^n$ when $c>-1$,  according to the literature. This allows to unify some formulas.

\subsection{The resolvents and the heat kernels of $B^d$ and $B^n$}

We start by recalling some well-known facts about the modified Bessel functions   $I_{\nu}$ and $K_{\nu}$
which constitute a basis of solutions of the modified Bessel equation
\begin{equation*}
\label{eq.mbessel}
z^2\frac{d^2v}{dz^2}+
z\frac{dv}{dz}-(z^2+\nu^2)v=0,
\quad \Rp z>0.
\end{equation*}
We recall that for $\Rp z>0$ one has 
\begin{equation*} \label{IK}
I_\nu(z)=\left (\frac{z}{2} \right )^\nu \sum_{m=0}^\infty\frac{1}{m!\,\Gamma (\nu+1+m)}\left (\frac{z}{2}\right )^{2m}, \quad K_\nu(z)=\frac{\pi}{2}\frac{I_{-\nu} (z)-I_\nu(z)}{\sin \pi \nu},
\end{equation*}
where limiting values are taken for the definition of $K_\nu$ when $\nu $ is an integer. The basic properties of these functions we need are collected in the following lemma, see e.g., \cite[Sections 9.6 and 9.7]{AS}.

\begin{lem}\label{behave}
For $\nu> -1$, $I_\nu$ is increasing and  $K_\nu$ is decreasing (when restricted to the positive real half line). Moreover they  satisfy
the following properties if $z \in \Sigma_{\pi/2-\eps}$.
\begin{itemize}
\item[(i)] $I_\nu(z)\neq 0$ for every  $\Rp z>0$. 
\item[(ii)] $I_{\nu}(z)\approx 
\frac{1}{\Gamma(\nu+1)} \left(\frac{z}{2}\right)^{\nu}, \quad \text{as }|z|\to 0,\qquad I_{\nu}(z)\approx \frac{e^z}{\sqrt{2\pi z}}(1+O(|z|^{-1}),\quad \text{as }|z|\to \infty$.\medskip
\item[(iii)] If $\nu\neq 0$,\quad  $K_{\mu}(z)\approx 
\frac{\nu}{|\nu|}\frac{1}{2}\Gamma(|\nu|) \left(\frac{z}{2}\right)^{-|\nu|}, \qquad K_{0}(z)\approx  -\log z, \qquad \text{as }|z|\to 0$
\\[2ex]
$ K_{\mu}(z)\approx \sqrt{\frac\pi{2z}}e^{-z},\quad \text{as }|z|\to \infty$.\medskip
\item[(iv)] $ I_{\nu}'(z)=I_{\nu+1}(z)+\frac{\nu}{z}I_{\nu}(z)$,\quad  $ K_{\nu}'(z)=K_{\nu+1}(z)+\frac{\nu}{z}K_{\nu}(z)$, for every  $\Rp z>0$.\\[1ex]
\medskip

\end{itemize}
\end{lem}

Note that 
\begin{equation}\label{Asymptotic I_nu}
|I_\nu(z)|\simeq C_{\nu,\epsilon} (1\wedge |z|)^{\nu+\frac 1 2}\frac{e^{Re z}}{\sqrt {|z|}},\qquad z\in \Sigma_{\frac \pi 2-\epsilon}
\end{equation}
for suitable constants $C_{\nu,\epsilon}>0$ which may be different in lower an in the upper estimate. 

Let us compute the resolvent operator of $B^n$. When we write $\sqrt z$  we mean the  square root of $z$ having positive real part.

\begin{prop} \label{risolventeB^n}
Let $c>-1$ and $\lambda\in \C \setminus (-\infty,0]$. Then, for every $f\in L^2_c$, 
$$(\lambda -B^n)^{-1}f=\int_0^\infty G^n(\lambda,y,\rho)f(\rho) \rho^c d\rho$$ with
\begin{equation}  \label{resolvent}
G^n(\lambda,y,\rho)=\begin{cases}
 y^{\frac{1-c}{2}} \rho^{\frac{1-c}{2}}\,I_{\frac{c-1}{2}}(\sqrt{\lambda}\,y)K_{{\frac{|1-c|}{2}}}(\sqrt{\lambda}\,\rho)\quad y\leq \rho\\[1.5ex]
 y^{\frac{1-c}{2}} \rho^{\frac{1-c}{2}}\,I_{\frac{c-1}{2}}(\sqrt{\lambda}\,\rho)K_{\frac{|1-c|}{2}}(\sqrt{\lambda}\,y)\quad y\geq \rho,
\end{cases}
\end{equation}

\end{prop}
{\sc Proof.} Let us first consider the case $\lambda=\omega^2$, $|\omega|=1$. By setting $u(y)= y^{\nu}v( \omega y)$, $\nu=(1-c)/2$, the homogeneous equation
$$D_{yy}u+\frac{c}{y}D_yu-\omega^2 u=0$$ transforms into the complex equation
\begin{equation*}
z^2\frac{d^2v}{dz^2}+
z\frac{dv}{dz}-(z^2+\nu^2)v=0,
\quad \Rp z>0.
\end{equation*}
Assume first that $-1 <c \le 1$ so that $ 0 \le \nu <1$. Then  
$u_1(z)=z^{\nu} I_{-\nu}( \omega z)$ and $u_2(z)=z^{\nu} K_{\nu}( \omega z)$ constitute a basis of solutions.
Since the Wronskian of
$K_{\nu}$, $I_{-\nu}$ is $1/r$,  see \cite[9.6 and 9.7]{AS}, that of 
$u_1$,$u_2$ is $r^{-c}$.
It follows that every solution of  
$$D_{yy}u+\frac{c}{y}D_yu-\omega^2u=f$$
is given by
\begin{equation} \label{defu}
u(y)=\int_{0}^{\infty}G^n(\omega^2,y,\rho)f\left(\rho\right) \rho^c d\rho+c_1y^{\nu} I_{-\nu}(\omega y)+c_2y^{\nu} K_{\nu}(\omega y),
\end{equation}
with $c_1,\ c_2\in\C$
 and 
\begin{equation*}  
G^n(\omega^2,y,\rho)=\begin{cases}
 y^{\nu} \rho^{\nu}\,I_{-\nu}(\omega y)K_{\nu}(\omega \rho)\quad y\leq \rho\\[1ex]
 y^{\nu} \rho^{\nu}\,I_{-\nu}(\omega \rho)K_{\nu}(\omega y)\quad y\geq \rho
\end{cases}
\end{equation*}
Next we use  Lemma \ref{behave} to show that
$$\sup _{y\in (0,+\infty)}\int_{0}^{\infty}|G^n(\omega,y,\rho)| \rho^c d\rho<+\infty.$$ 
Indeed, for $y\leq 1$, recalling that $|\omega|=1$, one has
\begin{align*}
 \int_{0}^{\infty}|G^n(\omega^2,y,\rho)|\rho^c d\rho=&\int_{0}^{y}y^{\nu}\rho^{\nu}\,|I_{-\nu}( \omega \rho)||K_{\nu}( \omega y)|\,\rho^cd\rho+\int_{y}^{\infty}y^{\nu}\rho^{\nu}\,|I_{-\nu}( \omega y)||K_{\nu}( \omega \rho)|\,\rho^cd\rho\\
\le& C\left (\int_0^y \rho^{c}d\rho+\int_y^1\rho^cd\rho+ \int_1^\infty  \rho^\frac{1+c}{2}(\sqrt{\rho})^{-1} e^{-\Rp {\omega }\rho} \right ) \le C
\end{align*}
and similarly for $y>1$.
By the symmetry of the kernel and  Young's inequality the integral operator $T$ defined by  $G^n(\omega^2, \cdot, \cdot)$ is therefore
 bounded in $L^2_c$.

Let $f\in C_c^\infty((0,\infty))$ with support in $(a,b)$ with $a>0$ and $u=(\omega^2-B^n)^{-1}f\in D(B^n)$. Then $u$ is given by \eqref{defu} with $c_1=0$, since $T$ is bounded in $L^2_c$, $K_{\nu}$ is exponentially decreasing and $I_{-\nu}$ is exponentially increasing near $\infty$. Since
\begin{align*}
u(y)=&\int_{0}^y y^{\nu} \rho^{\nu}K_{\nu}( \omega y)I_{\nu}( \omega \rho)f\left(\rho\right)\,\rho^c d\rho+\int_{y}^b y^{\nu} \rho^{\nu} K_{\nu}( \omega \rho)I_{-\nu}( \omega y)f\left(\rho\right)\,\rho^c d\rho+c_2y^{\nu}K_{\nu}( \omega y)
\end{align*}
we have for $y<a$
\begin{align*}
u(y)=&\int_{a}^by^{\nu} \rho^{\nu} K_{\nu}( \omega \rho)I_{-\nu}( \omega y)f\left(\rho\right)\,\rho^c d\rho+c_2y^{\nu}K_{\nu}( \omega y)=c_1y^{\nu} I_{-\nu}(\omega y)+c_2y^{\nu} K_{\nu}(\omega y)
\end{align*}
for some $c_1, c_2 \in\C$.
From Lemma \ref{behave} it follows that $v(y)=y^{\nu} I_{-{\nu}}( \omega y)$ satisfies the Neumann condition $\lim_{y \to 0} y^c v'(y) =0$ whereas $y^{\nu} K_{\nu}( \omega y)$ does not. Since $u \in D(B^n)$, by Lemma \ref{Neumann} $y^c u'(y) \to 0$ and hence $c_2=0$.
By density, $(\omega^2-B^n)^{-1}=T$, since both operators are bounded and coincide on compactly supported functions.

Finally let us compute the resolvent for a general $\lambda \not \in (-\infty,0]$.  If $M_s u(y)=u(sy)$, then
$M_{\sqrt{|\lambda|}}B^n M_{\sqrt{|\lambda|}^{-1}}=\frac{1}{|\lambda|} B^n$; setting $\lambda=|\lambda|\omega$ we get using the previous step
 \begin{align*}
(\lambda-B^n)^{-1}f &= |\lambda|^{-1}M_{\sqrt{|\lambda|}}(\omega-B^n)^{-1}M_{\sqrt{|\lambda|}^{-1}}f=\frac{1}{|\lambda|}\int_{0}^{\infty}G^n(\omega ,y\sqrt{|\lambda|},\rho)f\left(\frac{\rho}{\sqrt{|\lambda|}}\right) \rho^c d\rho\\&=|\lambda|^{\frac{c-1}{2}}\int_{0}^{\infty}G^n(\omega,y\sqrt{|\lambda|},\rho\sqrt{|\lambda|})f\left(\rho\right)\,\rho^cd\rho
\end{align*}
which gives \eqref{resolvent} when $-1<c \le 1$. When $c>1$, we use $I_{|\nu|}$, $K_{|\nu|}$ as a basis of solutions of Bessel equation and proceed  as before.
\qed

A similar proof gives the resolvent of $B^d$. We omit the details, see also \cite[Section 4.2]{met-negro-spina 5}.

\begin{prop} \label{risolventeB^d}
Let $c<1$ and  $\lambda\in \C \setminus (-\infty,0]$. Then, for every $f\in L^2_c$, 
$$(\lambda -B^d)^{-1}f=\int_0^\infty G^d(\lambda,y,\rho)f(\rho) \rho^c d\rho$$ with
\begin{equation}  \label{resolvent}
G^d(\lambda,y,\rho)=\begin{cases}
 y^\frac{1-c}{2}\rho^\frac{1-c}{2}\,I_{\frac{1-c}{2}}(\sqrt{\lambda}\,y)K_{\frac{1-c}{2}}(\sqrt{\lambda}\,\rho)\quad y\leq \rho\\[2ex]
 y^\frac{1-c}{2}\rho^\frac{1-c}{2}\,I_{\frac{1-c}{2}}(\sqrt{\lambda}\,\rho)K_{\frac{1-c}{2}}(\sqrt{\lambda}\,y)\quad y\geq \rho.
\end{cases}
\end{equation}
\end{prop}

Note that  when $|c|<1$ the resolvent of $B^n$ uses $I_{\frac{c-1}{2}}$ whereas $B^d$ is constructed with  $I_{\frac{1-c}{2}}$.

\smallskip

Next we compute the heat kernel of $e^{z B^n}, e^{z B^d}$, proceeding as in \cite[Section 4.2]{met-negro-spina 5}. These heat kernels are known and usually computed by probabilistic methods. Instead we provide a purely analytical proof and refer also to \cite{Garofalo} for a similar approach in the Neumann case.



For $z\in C_+, y,\rho>0$ we denote now by $p(z,y,\rho)$ the heat kernel of the operator $B$ and argue first for positive $t$. We look for a smooth function $p(t,y,\rho)$ such that, for every $f \in L^2_c$ 
$$
e^{tB}f(y)=\int_0^\infty p(t,y,\rho)f(\rho)\, d\rho.
$$
Note that the kernel is written with respect to the Lebesgue measure rather than $y^c dy$.  The function $p$ should then satisfy
\begin{equation*} 
\begin{cases}
 p_t(t,y,\rho)=D_{yy}p(t,y,\rho)+\frac{c}{y}D_yp(t,y,\rho)\\
 p(0,y,\rho)=\delta_\rho.
\end{cases}
\end{equation*}
It follows that $\tilde{p}(t,y,\rho)=y^\frac{c}{2}p(t,y,\rho)\rho^{-\frac{c}{2}}$  satisfies with $\nu^2=(c-1)^2/4$
\begin{equation} \label{parabolic}
\begin{cases}
 {\tilde p}_t(t,y,\rho)=D_{yy}{\tilde p}(t,y,\rho)-\frac{1}{y^2}\left(\nu^2-\frac{1}{4}\right){\tilde p}(t,y,\rho)\\
 {\tilde p}(0,y,\rho)=\delta_\rho.
\end{cases}
\end{equation}
Since $\lambda^2B=M_\lambda^{-1}B M_\lambda$ we obtain $e^{t\lambda^2B}=M_\lambda^{-1}e^{tB}M_\lambda$. Rewriting this identity using the kernel $\tilde p$ and setting $\lambda^2 t=1$ we obtain
$${\tilde p}(t,y,\rho)=\frac{1}{\sqrt{t}}{\tilde p}\left(1,\frac{y}{\sqrt{t}},\frac{\rho}{\sqrt{t}}\right):=\frac{1}{\sqrt{t}}F\left(\frac{y}{\sqrt{t}},\frac{\rho}{\sqrt{t}}\right).$$
Then  (\ref{parabolic}) becomes 
\begin{align*}
D_{yy}F\left(\frac{y}{\sqrt{t}},\frac{\rho}{\sqrt{t}}\right)&-\frac{1}{y^2}\left(\nu^2-\frac{1}{4}\right)tF\left(\frac{y}{\sqrt{t}},\frac{\rho}{\sqrt{t}}\right)+\\&+\frac{1}{2}F\left(\frac{y}{\sqrt{t}},\frac{\rho}{\sqrt{t}}\right)+\frac{1}{2}\frac{y}{\sqrt{t}}D_yF\left(\frac{y}{\sqrt{t}},\frac{\rho}{\sqrt{t}}\right)+\frac{1}{2}\frac{\rho}{\sqrt{t}}D_{\rho}F\left(\frac{y}{\sqrt{t}},\frac{\rho}{\sqrt{t}}\right)=0
\end{align*}
that is
\begin{align*}
D_{yy}F\left(y,\rho\right)&-\frac{1}{y^2}\left(\nu^2-\frac{1}{4}\right)F\left(y,\rho\right)+\frac{1}{2}F\left(y,\rho\right)+\frac{1}{2}yD_yF\left(y,\rho\right)+\frac{1}{2}\rho D_{\rho}F\left(\frac{y}{\sqrt{t}},\frac{\rho}{\sqrt{t}}\right)=0.
\end{align*}
Since for large $y$ the operator $B$ behaves like $D^2$, 
having in mind the gaussian kernel, we look for  a solution of the form
$$F(y,\rho)=\frac{1}{\sqrt {4 \pi}}\exp\left\{-\frac{(y-\rho)^2}{4}\right\}H(y\rho)$$
with $H$ depending only on the product of the variables.
By straightforward computations, we deduce
$$\rho^2D_{yy}H(y\rho)+\rho^2D_yH(y\rho)-\frac{1}{y^2}\left(\nu^2-\frac{1}{4}\right)H(y\rho)=0$$
or
$$DH(x)+D_{x}H(x)-\frac{1}{x^2}\left(\nu^2-\frac{1}{4}\right)H(x)=0.$$
Setting $H(x)=u(x)e^{-\frac{x}{2}}$, $u$ solves
$$D_{xx}u-\frac{1}{4}u(x)-\frac{1}{x^2}\left(\nu^2-\frac{1}{4}\right)u(x)=0$$ and $v(x)=u(2x)$ satisfies
$$D_{xx}v-v(x)-\frac{1}{x^2}\left(\nu^2-\frac{1}{4}\right)v(x)=0.$$
It follows that $v(x)=c_1\sqrt{x}I_{\nu}(x)+c_2\sqrt{x} K_{|\nu|} (x)$. Since the function $H$ captures the behaviour of the heat kernel near the origin (the behaviour at infinity is governed by the gaussian factor) and since the resolvents of $B^n, B^d$ are constructed with $\nu=(c-1)/2$, $\nu=(1-c)/2$, respectively,   we choose $ c_2=0$, $c_1=1$  and $\nu$ accordingly. Therefore  in the case of Neumann boundary conditions,
$u(x)=v\left(\frac{x}{2}\right)=\kappa\sqrt{\frac{x}{2}}I_{\frac{c-1}{2}}\left(\frac{x}{2}\right)$, $H(y\rho)=u(y\rho)e^{-\frac{y\rho}{2}}=\kappa\sqrt{\frac{y\rho}{2}}I_{\frac{c-1}{2}}\left(\frac{y\rho}{2}\right)e^{-\frac{y\rho}{2}}$, 
$$F(y,\rho)=\frac{\kappa}{\sqrt{4\pi}}\exp\left\{-\frac{(y-\rho)^2}{4}\right\}\sqrt{\frac{y\rho}{2}}I_{\frac{c-1}{2}}\left(\frac{y\rho}{2}\right)e^{-\frac{y\rho}{2}}=\frac{\kappa}{\sqrt{4\pi}}\sqrt{\frac{y\rho}{2}}\exp\left\{-\frac{y^2+\rho^2}{4}\right\}I_{\frac{c-1}{2}}\left(\frac{y\rho}{2}\right)$$ and
\begin{equation*} 
\tilde{p}(t,y,\rho)=\frac{1}{\sqrt{4\pi t}}H\left (\frac{y\rho}{t}\right)\exp\left\{-\frac{(y-\rho)^2}{4t}\right\}=\frac{\kappa}{t\sqrt{4\pi}}\sqrt{\frac{y\rho}{2}}\exp\left\{-\frac{y^2+\rho^2}{4t}\right\}I_{\frac{c-1}{2}}\left(\frac{y\rho}{2t}\right).
\end{equation*}
It follows that
\begin{equation} \label{defp}
p^n(t,y,\rho)=y^{-\frac{c}{2}}\tilde{p}(t,y,\rho)\rho^{\frac{c}{2}}=\frac{\kappa}{t\sqrt{4\pi}}\left(y\rho\right)^{\frac{1-c}{2}}\rho^c\exp\left\{-\frac{y^2+\rho^2}{4t}\right\}I_{\frac{c-1}{2}}\left(\frac{y\rho}{2t}\right).\end{equation}
In the case of Dirichlet boundary conditions it is sufficient to change $I_{\frac{c-1}{2}}$ with $I_{\frac{1-c}{2}}$ to obtain the corresponding kernel.
\begin{equation} \label{defpDirichlet}
p^d(t,y,\rho)=y^{-\frac{c}{2}}\tilde{p}(t,y,\rho)\rho^{\frac{c}{2}}=\frac{\kappa}{t\sqrt{4\pi}}\left(y\rho\right)^{\frac{1-c}{2}}\rho^c\exp\left\{-\frac{y^2+\rho^2}{4t}\right\}I_{\frac{1-c}{2}}\left(\frac{y\rho}{2t}\right).\end{equation}

Finally, we give a formal proof and compute the constant $\kappa$.

\begin{teo} \label{preciseKernels} For $z \in \C_+$ the heat kernels of the operators $B^n, B^d$ are given by 
$$p_{B^n}(z,y,\rho)=\frac{1}{2z}\rho^c(y\rho)^\frac{1-c}{2}\exp\left\{-\frac{y^2+\rho^2}{4z}\right\}I_{\frac{c-1}{2}}\left(\frac{y\rho}{2z}\right),\qquad c>-1.$$
$$
 p_{B^d}(z, y,\rho)=\frac{1}{2z}\rho^{c}\left (y\rho \right)^\frac{1-c}{2}\exp\left\{-\frac{y^2+\rho^2}{4z}\right\}I_{\frac{1-c}{2}}\left(\frac{y\rho}{2z}\right),\qquad c<1.
$$
\end{teo}

{\sc Proof.} Let us consider  $B^n$ and $t\in \R^+$, first. The Laplace transform of the right hand side of \eqref{defp} is given by, see \cite[p.200]{Erdelyi},
\begin{equation*}
\begin{cases}
\frac {2\kappa}{\sqrt{4\pi}}\rho^c(y\rho)^\frac{1-c}{2} I_{-\sqrt{D}}(y\sqrt{\lambda})K_{\sqrt{D}}(\rho\sqrt{\lambda}) \quad y\leq \rho \\[1ex]
\frac{2\kappa}{\sqrt{4\pi}}\rho^c(y\rho)^\frac{1-c}{2} I_{-\sqrt{D}}(\rho\sqrt{\lambda})K_{\sqrt{D}}(y\sqrt{\lambda}) \quad y\geq \rho. 
\end{cases}
\end{equation*}
For $\kappa=\sqrt{\pi}$ it coincides with the kernel $\rho^cG^n(\lambda, y,\rho)$ of the resolvent operator $(\lambda-B^n)^{-1}$, see Proposition \ref{risolventeB^n}. Let $S(t)$ be the operator defined through the kernel $p_{B^n}$, that is \eqref{defp}  with $\kappa=\sqrt{\pi}$. By  Lemma \ref{boundedlpm theta} below,  $\|S(t)\| \le C$, $t \ge 0$,  in $L^2_c$.
Given $f \in C_c^\infty ((0,\infty))$, let $u(t,y)=S(t)f(y)$. By the construction of the kernel $p$ we have $u_t=B u$ pointwise. Finally, for $\lambda>0$,
\begin{align*}
\int_0^\infty e^{-\lambda t}u(t,y)\, dt&=\int_0^\infty e^{-\lambda t}\, dt \int_0^\infty p_{B^n}(t,y,\rho)f(\rho)\, d\rho=\int_0^\infty f(\rho)\, d\rho \int_0^\infty e^{-\lambda t}p_{B^n}(t,y,\rho)\, dt \\
&=\int_0^\infty G^n(\lambda,y,\rho)\rho^cf(\rho)\, d\rho.
\end{align*}
It follows that the Laplace transform of $S(t)f$ coincides with the resolvent of $B^n$, hence, by uniqueness, $S(t)$ is the generated semigroup and $p_{B^n}(t, \cdot, \cdot)$ its kernel.

For complex times we argue in a similar way;   we fix $0\leq|\theta|<\frac \pi 2$ and $\omega=e^{i\theta}\in \C_+$. Then for $t>0$, $p(t\omega,\cdot,\cdot)$ is the heat kernel of the scaled semigroup $T_\omega(t)=e^{t\omega B^n}$ whose generator is $A_\omega=\omega B^n$. Its resolvent  is then given, for $\lambda>0$,  by $(\lambda -A_\omega)^{-1}=\omega^{-1} \left(\omega^{-1}\lambda-B^n\right)^{-1}$ and its integral kernel is $\omega^{-1}G^n(\omega^{-1}\lambda,y\rho)$. The same argument  above applied to $T_\omega(t)$ proves then  the assertion for $z=t\omega$.
The proof for  $B^d$ is similar.
\qed
 The following result is proved in \cite{Garofalo1}.
\begin{prop}
If $c>-1$, then $e^{zB^n}1=1$.
\end{prop}

{\sc{Proof.}} The proof follows  using the explicit expression of $p_{B^n}$ of Theorem \ref{preciseKernels} and  the identity 
\begin{align*}
\int_0^\infty e^{-\alpha\rho^2}\rho^{1+\nu}I_\nu(z\rho)\,d\rho=\frac{z^\nu}{(2\alpha)^{\nu+1}}\,e^{\frac{z^2}{4\alpha}},\qquad \nu>-1
\end{align*}
which holds for every $z,\alpha\in \C_+$. See  \cite[Formula 11.4.29, page 486]{AS} where the latter equality is expressed in terms of the Bessel functions $J_\nu$ which satisfies $I_\nu(z)=e^{-\frac 1 2 \nu\pi i} J_\nu\left(e^{\frac 1 2 \pi i}z\right)$.\qed

\subsection{Heat kernel bounds}
The asymptotic behaviour of Bessel functions allows to deduce explicit bounds for the heat kernels $p_{B^n}$ and $p_{B^d}$. 
\begin{prop}\label{Estimates Bessel kernels}
Let $p_{B^n}, p_{B^d}$ be the kernels defined in Theorem \ref{preciseKernels}. Then for every $\eps>0$, there exist $C_\eps, \kappa_\eps>0$  such that 
for $z\in\Sigma_{\frac{\pi}{2}-\eps}$
\begin{align*} 
| p_{B^d}(z,y,\rho)|
&\leq 
 C_\eps |z|^{-\frac{1}{2}} \left (\frac{y}{|z|^{\frac{1}{2}}}\wedge 1 \right)^{1-c} \left (\frac{\rho}{|z|^{\frac{1}{2}}}\wedge 1 \right)
\exp\left(-\frac{|y-\rho|^2}{\kappa_\eps |z|}\right),  \quad c <1 
\end{align*}
and 
\begin{align*} 
| p_{B^n}(z,y,\rho)|
&\leq 
C_\eps |z|^{-\frac{1}{2}}  \left (\frac{\rho}{|z|^{\frac{1}{2}}}\wedge 1 \right)^{c}
\exp\left(-\frac{|y-\rho|^2}{\kappa_\eps |z|}\right), \quad c >-1 
\end{align*}
\end{prop}
{\sc{Proof.}} Using   \eqref{Asymptotic I_nu} we get 
\begin{align*} 
| p_{B^d}(z,y,\rho)|& \le  C_\eps\,  |z|^{-1}\rho^{c}\left (y\rho \right)^{\frac{1-c}{2}}\exp\left\{-\frac{Re z}{4 |z|^2}(y^2+\rho^2)\right\}\left(\frac{y\rho}{2|z|} \wedge 1\right)^{\frac{1-c}{2}+\frac 1 2}\left(\frac{2|z|}{y\rho}\right)^{\frac 1 2}\exp\left\{\frac{Rez}{2|z|^2}y\rho \right\}.\\
 &\le  C_\eps\, |z|^{-1}  \rho^{c}\left (y\rho \right)^{\frac{1-c}2}
\left ( \frac{y\rho}{|z|} \wedge 1 \right)^{1-\frac c2}
\left(\frac{2|z|}{y\rho}\right)^{\frac 1 2}\exp\left(-\frac{|y-\rho|^2}{\kappa'_\eps |z|}\right)\\[1ex]
\nonumber&=
C_\eps |z|^{-\frac 1 2 } \left (\frac y\rho \right)^{-\frac{c}2} \left ( \frac{y\rho}{|z|} \wedge 1 \right)^{1-\frac c 2}
\exp\left(-\frac{|y-\rho|^2}{\kappa'_\eps |z|}\right) \\[1ex]
\nonumber &\le C'_\eps |z|^{-\frac 1 2 }  \left ( \frac{y}{|z|^{-\frac 1 2 }} \wedge 1\right)^{1-c}\left ( \frac{\rho}{|z|^{-\frac 1 2 }} \wedge 1 \right)
\exp\left(-\frac{|y-\rho|^2}{\kappa_\eps |z|}\right) 
\end{align*}
by Lemmas \ref{equiv}, \ref{equiv1}. The proof for $p_{B^n}$ is similar.
\qed

Note that the constant $\kappa_\eps$ above is explicit. For example, for $z \ge 0$, that is $\eps=\pi/2$, then $\kappa'_\eps=4$ in the above proof and we can take $\kappa_\eps=4+\delta$ for every $\delta>0$.

Next we prove bounds for the gradients of the heat kernels.

\begin{prop}\label{Estimates gradient kernel bessel}
For every $\eps>0$, there exist $C_\eps, \kappa_\eps>0$  such that
for $z\in\Sigma_{\frac{\pi}{2}-\eps}$ 
\begin{align*} 
| D_y p_{B^d}(z,y,\rho)|
\leq 
C_\eps |z|^{-1} \left (\frac{y}{|z|^{\frac{1}{2}}}\wedge 1 \right)^{-c} \left (\frac{\rho}{|z|^{\frac{1}{2}}}\wedge 1 \right)
\exp\left(-\frac{|y-\rho|^2}{\kappa_\eps |z|}\right),\quad  c<1;
\end{align*}
and 
\begin{align*} 
| D_y p_{B^n}(z,y,\rho)|
&\leq 
C_\eps |z|^{-1}\left (\frac{y}{|z|^{\frac{1}{2}}}\wedge 1 \right)\left (\frac{\rho}{|z|^{\frac{1}{2}}}\wedge 1 \right)^{c}
\exp\left(-\frac{|y-\rho|^2}{\kappa_\eps |z|}\right), \quad c>-1.
\end{align*}
\end{prop}
{\sc{Proof.}} We give a proof first for $B^d$. Differentiating $p_{B^d}$ with respect to $y$ we obtain
\begin{align*}
D_y p_{B^d}(z,y,\rho)=\Bigg[\frac {1-c}{2 y}-\frac{y}{2z}+\frac{\rho}{2z}\frac{I'_{\nu}\left(\frac{y\rho}{2z}\right)}{I_{\nu}\left(\frac{y\rho}{2z}\right)}\Bigg]p_{B^d}(z,y,\rho),
\end{align*}
where   $\nu=(1-c)/2$.
We recall now that, see Lemma \ref{behave}, 
\begin{align*}
 I_{\nu}'(s)=I_{\nu+1}(s)+\frac{\nu}{s}I_{\nu}(s).
 \end{align*}
This implies 
\begin{align*}
D_y p_{B^d}(z,y,\rho)&=\Bigg[\frac {1-c}{ y}-\frac{y}{2z}+\frac{\rho}{2z}\frac{I_{\nu+1}\left(\frac{y\rho}{2z}\right)}{I_{\nu}\left(\frac{y\rho}{2z}\right)}\Bigg]p_{B^d}(z,y,\rho) \\
&=\frac{1}{\sqrt z}\left [(1-c)\frac {\sqrt z}{ y}+\frac{y}{2\sqrt z}\left(\frac{I_{\nu+1}\left(\frac{y\rho}{2z}\right)}{I_{\nu}\left(\frac{y\rho}{2z}\right)}-1\right)-\frac{y-\rho}{2\sqrt z}\frac{I_{\nu+1}\left(\frac{y\rho}{2z}\right)}{I_{\nu}\left(\frac{y\rho}{2z}\right)}\right ]p_{B^d}(z,y,\rho)
\end{align*}
The boundedness of $\frac{I_{\nu+1}}{I_{\nu}}$ gives 
$$\left |\frac{y-\rho}{2\sqrt z}\frac{I_{\nu+1}\left(\frac{y\rho}{2z}\right)}{I_{\nu}\left(\frac{y\rho}{2z}\right)}\right | \le C  \left |\frac{y-\rho}{2\sqrt z} \right | \le Ce^{\eps \frac{|y-\rho|^2}{|z|}}
$$
Next we use the estimate 
$\left |\frac{I_{\nu+1}(w)}{I_{\nu}(w)}-1 \right | \le C(1\wedge |w|^{-1})$ for
 $w\in\Sigma_{\frac{\pi}{2}-\eps}$ 
to bound 
$$K(\xi, \eta)=\xi \left(\frac{I_{\nu+1}\left(\frac{\xi \eta}{2}\right)}{I_{\nu}\left(\frac{\xi\eta}{2}\right)}-1\right),$$
where $\xi=\frac{y}{\sqrt z}$ and $\eta=\frac{\rho}{\sqrt z}$.

Clearly $|K(\xi, \eta)| \le C$ if $|\xi| \le 2$ and $|K(\xi, \eta)| \le Ce^{\eps |\xi-\eta|^2}$ if $|\xi| \ge 2$ and $|\eta | \le 1$. Finally, if $|\xi| \ge 2$, $|\eta | \ge 1$, then $|K(\xi, \eta)| \le C\frac{|\xi|}{|\xi \eta|} \le \frac{C}{|\eta |} \le C$. Then 
\begin{equation}  \label{Grad-etaGrande}
|D_y p_{B^d}(z,y \rho)| \le C\frac{1}{\sqrt{|z|}} \left (1+\frac{|(1-c)\sqrt z|}{y}\right )e^{\eps\frac{|y-\rho|^2}{|z|}}|p_{B^d}(z,y,\rho)|
\end{equation}
 and the thesis  follows from Proposition \ref{Estimates Bessel kernels}. 

Concerning $B^n$ we first note that  $\nu=(c-1)/2$  in the above notation.  Then we get
$$D_y p_{B^n}(z,y,\rho)=\Bigg[-\frac{y}{2z}+\frac{\rho}{2z}\frac{I_{\nu+1}\left(\frac{y\rho}{2z}\right)}{I_{\nu}\left(\frac{y\rho}{2z}\right)}\Bigg]p_{B^n}(z,y,\rho)=\frac{y}{2z}\Bigg[-1+\frac{\rho}{y}\frac{I_{\nu+1}\left(\frac{y\rho}{2z}\right)}{I_{\nu}\left(\frac{y\rho}{2z}\right)}\Bigg]p_{B^n}(z,y,\rho).$$ 
Proceeding as before we get \eqref{Grad-etaGrande} without the term $(1-c)\sqrt{z}/y$ and we only  look for a  better estimate in the region $\frac{y}{\sqrt{|z|}}<1$.
Setting $\xi=\frac{y}{\sqrt z}$ and $\eta=\frac{\rho}{\sqrt z}$ as before, we prove a bound for
$$K_0(\xi,\eta)=-1+\frac{\eta}{\xi}\frac{I_{\nu+1}\left(\frac{\eta\xi}{2}\right)}{I_{\nu}\left(\frac{\eta\xi}{2}\right)}$$
in the case $|\xi|<1$. 
Using the estimate $\frac{|I_{\nu+1}(w)|}{|I_{\nu}(w)|} \le C(1\wedge |w|)$, we get $|K_0(\xi,\eta)|\leq 1+C\left|\frac{\eta}{\xi}\right|\left(1\wedge |\eta\xi|\right)$.
Assume first $|\xi\eta|\leq 1$. Then
$|K_0(\xi,\eta)|\leq   C\left(1+|\eta|^2\right)\leq C$ if $|\eta|\leq 1$ and $|K_0(\xi, \eta)| \le Ce^{\eps |\xi-\eta|^2}$ if $|\eta| >1$.
Let now $|\xi\eta|> 1$.  Then $\left|\frac{\eta}{\xi}\right|\leq |\eta|^2$ and 
$|K_0(\xi,\eta)|\leq   C\left(1+|\eta|^2\right)\leq Ce^{\eps |\xi-\eta|^2}$.
It follows that, when $\frac{y}{\sqrt{|z|}}<1$,
$$|D_y p_{B^n}(z,y \rho)| \le \frac{C}{\sqrt{|z|}}\frac{y}{\sqrt{|z|}}e^{\eps |\xi-\eta|^2} |p_{B^n}(z,y,\rho)|$$  
and the thesis follows from Proposition \ref{Estimates Bessel kernels}. 
\qed

\section{The semigroups $e^{zB^d}$ and $e^{zB^n}$}
In this section we show that the operators defined  through the kernels $p_{B^d}$  and  $p_{B^n}$ are  strongly continuous semigroups in $L^p_m=L^p(\R_+; y^mdy)$.
We define $\{e^{zB}\}_{z\in\C_+}$,  $\{D_y e^{zB}\}_{z\in\C_+}$  for $f\in C_c^\infty(0,\infty)$  by 
\[
 [e^{zB}f](y):=
\int_0^\infty
  p(z,y,\rho)f(\rho)
\,d\rho, \quad  [D_y e^{zB}f](y):=
\int_0^\infty
D_y  p(z,y,\rho)f(\rho)
\,d\rho
\]
where $p=p_{B^d}$ or $p=p_{B^n}$ and, accordingly, we write $e^{z B^d}$ and $e^{z B^n}$.

The following lemma is consequence of the heat kernel estimates of Propositions \ref{Estimates Bessel kernels}, \ref{Estimates gradient kernel bessel} and Proposition \ref{Boundedness theta}.

\begin{lem} \label{boundedlpm theta} 
Let  $\theta\geq 0$, $\delta =\pi/2-\eps$, $\eps>0$. The following properties hold for $z \in \Sigma_\delta$.\smallskip
\begin{itemize}
\item[(i)] If $c<1$ and $c-1+\theta<\frac{m+1}{p}<2$, then
$\|e^{zB^d}\|_{L^p_m\to L^p_{m-p{\theta}}}\leq C|z|^{-\frac{\theta} 2}.$\smallskip
\item[(ii)] If $c<1$ and  $c+\theta<\frac{m+1}{p} <2$, then
$
\|\sqrt {z }D_ye^{zB^d}\|_{L^p_m\to L^p_{m-p{\theta}}}\leq C|z|^{-\frac{\theta} 2}.
$\smallskip
\item[(iii)] If $c>-1$ and $\theta<\frac{m+1}{p} <c+1$, then  
$\|e^{zB^n}\|_{L^p_m\to L^p_{m-p{\theta}}}\leq C|z|^{-\frac{\theta} 2}.$\smallskip
\item[(iv)] If $c>-1$ and $\theta-1<\frac{m+1}{p} <c+1$, then  
$\|\sqrt z D_ye^{zB^n}\|_{L^p_m\to L^p_{m-p{\theta}}}\leq C|z|^{-\frac{\theta} 2}.$
\end{itemize}
\end{lem}


\begin{prop}\label{Gen DN}
\begin{itemize}
\item[(i)] If $c<1$ and $c-1< \frac{m+1}{p} <2$, then  $\{e^{zB^d}\}$ is a bounded analytic semigroup of angle $\pi/2$ in $L^p_m$.
\item[(ii)]If $c>-1$ and $0< \frac{m+1}{p} <c+1$, then  $\{e^{zB^n}\}$ is a bounded analytic semigroup of angle $\pi/2$ in $L^p_m$.
\end{itemize}
\end{prop}
{\sc Proof. } The boundedness of the families $\{e^{zB^d}\}_{ z\in \Sigma_\delta}$, $\{e^{zB^n}\}_{ z\in \Sigma_\delta}$ follows from the previous lemma with $\theta=0$; the semigroup law is inherited  from the one of $L^2_c$ via a density argument and we have only to prove the  strong continuity at $0$. Let $f, g \in C_c^\infty (0, \infty)$. Then as $z \to 0$, $z \in \Sigma_\delta$, 
$$
\int_0^\infty (e^{zB}f)\, g\, y^m dy=\int_0^\infty (e^{zB}f) \,g\, y^{m-c}y^c  dy \to \int_0^\infty fgy^{m-c}y^cdy  =\int_0^\infty fgy^{m}dy,
$$
by the strong continuity of $e^{zB}$ in $L^2_c$. By density and uniform boundedness of the family $(e^{zB})_{ z\in \Sigma_\delta}$ this holds for every $f \in L^p_m$, $g \in L^{p'}_m$. The semigroup is then weakly continuous, hence strongly continuous.\qed

We denote by $B_{m,p}^d, B_{m,p}^n$ the generators of $e^{zB^d}, e^{zB^n}$ in $L^p_m$ and characterize their domain. Observe that, since the  heat kernels of these semigroups are given by Theorem \ref{preciseKernels}, their resolvents are those of Propositions \ref{risolventeB^d}, \ref{risolventeB^n}. \\
We recall that the Sobolev spaces $W^{k,p}_m$ are studied in detail in Appendix B and that traces at the boundary $y=0$ are well-defined when $(m+1)/p<1$.

It is useful to define $D(B_{m,p, max})=\{u \in L^p_m \cap W^{2,p}_{loc}(\R_+): Bu \in L^p_m\}$. We start with $B^n$.

\begin{prop} \label{dominioN}
If $c>-1$ and $0<\frac{m+1}{p}<c+1$, then 
\[
D(B_{m,p}^n)=
\{
u\in D(B_{m,p, max}): \;
\frac{D_y u}{y}, \;
D_{yy}u\in L^p_m 
\}.
\]
Moreover, $D(B_{m,p}^n)=\{u \in W^{2,p}_m: D_y u (0)=0 \}$ when $\frac{m+1}{p} <1$ and  $D(B_{m,p}^n)= W^{2,p}_m$ when $\frac{m+1}{p} >1$.
\end{prop}
{\sc Proof. } Let $D$ be the right-hand side above,  $u\in D(B^d_{m,p})$ and let $f:=u-B^n_{m,p}u$.  Then $u=(I-B^n_{m,p})^{-1}f=\int_0^\infty e^{-t}e^{tB^n}f\,dt$ and  $D_yu=\int_0^\infty e^{-t}D_ye^{tB^d}f\,dt$.  Using Minkowski's inequality  and  Lemma \ref{boundedlpm theta} (iv) with  $\theta-1<0<\frac{m+1}{p} <c+1$,  we get 
\begin{align*}
\|y^{-\theta}D_yu\|_{L^p_m}=\|D_yu\|_{L^p_{m-\theta p}}\leq \int_0^\infty e^{-t}\frac{1}{\sqrt t}\|\sqrt tD_ye^{tB^d}f\|_{L^p_{m-\theta p}}\,dt\leq C\|f\|_{L^p_{m}}\int_0^\infty e^{-t}t^{-\frac{\theta+1}{2}}\,dt\,
 \end{align*}
and then $ y^{-\theta}D_y u \in L^p_m$ for every $0 \le \theta<1$. To reach $\theta=1$ let $v=D_y u$ and $g=D_y v+c\frac{v}{y}=Bu$. Then 
\begin{equation} \label{derivata}
v(y)=y^{-c}\int_0^y g(s)s^c\, ds +Ky^{-c}:=w(y)+Ky^{-c}
\end{equation}
(note that the integral converges, by H\"older inequality, since $(c+1)>\frac{m+1}{p}$). By Hardy inequality, see Lemma \ref{Hardy}, $\|\frac{w}{y}\|_{L^p_m} \le C\|g\|_{L^p_m}$ but $y^{-c-\theta} \not \in L^p_m(0,1)$ for $\theta<1$, sufficiently close to $1$. It follows that $K=0$, $v=w$ and then $D_y u /y \in L^p_m$ and, by difference, $D_{yy}u \in L^p_m$, too.

This shows that $D(B^n_{m,p}) \subset D$ and we have only to show that $I-B$ is injective on $D$. Assume that $u \in D$ and that $u-Bu=0$, then $u(y)=c_1y^{\frac{1-c}{2}}I_{\frac{c-1}{2}}+c_2 y^{\frac{1-c}{2}}K_{\frac{|1-c|}{2}}$. However $c_1=0$, since $I_{\frac{c-1}{2}}$ is exponentially increasing at $\infty$. Concerning $u_2(y)= y^{\frac{1-c}{2}}K_{\frac{|1-c|}{2}}$ we note that its derivative, as $y \to 0$,  behaves like $y^{-1}$ when $c\le 1$ and like $y^{-c}$ when $c>1$. In both cases $(D_y u_2)/y \not \in L^p_m$. Then $c_2=0$ and $u=0$.

The last part follows from Proposition \ref{notrace} applied to $D_y u$, taking into account that $W^{1,p}_m=W^{1,p}_{0,m}$ when $(m+1)/p>1$. In the case $(m+1)/p<1$ observe that $D_y u$ has a finite limit as $y \to 0$, by Lemma \ref{l1},  which must be equal to $0$, otherwise $(D_y u)/y \not \in L^p_m$.
\qed

\begin{cor} \label{dominioN1} If $c>-1$,  $0<\frac{m+1}{p}<c+1$  and  $u \in D(B^n_{m,p})$, then  $\lim_{y \to 0}y^{\frac{m+1}{p}-1}D_y u=0$, hence $\lim_{y \to 0}y^c D_y u=0$. Moreover
\begin{itemize}
\item[(i)] if $\frac{m+1}{p}<2$, then $D_y u \in L^1(0,1)$ and therefore $ \lim_{y \to 0} u(y)$ exists finite;
\item[(ii)] if $\frac{m+1}{p}=2$, then $\frac{u}{y^{2\theta}} \in L^p_m$ if $2\theta <2$;
\item[(iii)] if $\frac{m+1}{p}>2$ then $\frac{u}{y^2} \in L^p_m$.
\end{itemize}
\end{cor}
{\sc Proof. }  By \eqref{derivata} with $K=0$ and  H\"older inequality we get
$$|D_y u|\le y^{-c} \int_0^y |g(s)|s^{\frac mp}s^{c-\frac mp}\, ds \le y^{1-\frac{m+1}{p}}\left (\int_0^y |g(s)|^p s^m\, ds \right )^{\frac1p}
$$ and the first statement follows and yields $D_y u \in L^1(0,1)$ when $(m+1)/p<2$. This proves $(i)$. The proof of $(ii)$ is similar since $D_yu = o(y^{-1})$, hence $u =o(\log y)$  as $y \to 0$.

Assume now that $(iii)$ holds and write, by \eqref{derivata} again,
$$D_yu=y^{-c} \int_0^y g(s)s^c\, ds=y\int_0^1 g(ty)t^c\, dt.$$ Then
$$
\frac{|u(y)-u(1)|}{y^2} \le \frac{1}{y^2} \int_y^1 s\, ds \int_0^1 |g(ts)|t^c\, dt \le \int_1^\infty \eta\, d\eta \int_0^1 |g(t\eta y)| t^c\, dt
$$
and Minkowski's inequality gives
$$
\left \|\frac{|u(y)-u(1)|}{y^2} \right \|_{L^p_m} \le \int_1^\infty \eta\, d\eta \int_0^1 t^c \|g(t\eta \cdot)\|_{L^p_m}\, dt=\|g\|_{L^p_m} \int_1^\infty \eta^{1-\frac{m+1}{p}} d \eta \int_0^1 t^{c-\frac{m+1}{p}} dt.
$$
Since also $y^{-2} u(1) \in L^p_m (0,1)$, the proof of $(iii)$ is complete.
\qed

Next we consider $B^d$. 

\begin{prop}\label{Dirichlet}
 Let $c<1$ and $c-1<\frac{m+1}p<2$. Then 
\begin{eqnarray*}
D(B_{m,p}^d)=&
 \{
u\in D(B_{m,p, max}): \;
 y^{-2\theta}u\in L^p_m {\rm\  for\  every\ } 0\leq\theta \le 1 \\[1ex]
&\ {\rm  such\  that\ } c-1+2\theta<\frac{m+1}{p} <2 \}. 
\end{eqnarray*}
Moreover
\begin{itemize}
\item[(i)] $(1\wedge u)^{2-2\theta}D_{yy}u, (1\wedge u)^{1-2\theta}D_y u \in L^p_m$ for every $u \in D(B_{m,p}^d)$ and $\theta$ as above;
\item[(ii)]
$D_y u \in L^1(0,1)$ and $\lim_{y \to 0}u(y)=0$ for every $u \in D(B_{m,p}^d)$.
\end{itemize}

\end{prop}
{\sc{Proof.}}  Let $u\in D(B^d_{m,p})$ and let $f:=u-B^d_{m,p}u$.  Then $u=(I-B^d_{m,p})^{-1}f=\int_0^\infty e^{-t}e^{tB^d}f\,dt$. Then using Minkowski's inequality and  Lemma \ref{boundedlpm theta} we get when  $ c-1+2\theta<\frac{m+1}{p} <2 $ and $0\leq\theta<1$, 
\begin{align*}
\|y^{-2\theta}u\|_{L^p_m}=\|u\|_{L^p_{m-2\theta p}}\leq \int_0^\infty e^{-t}\|e^{tB^d}f\|_{L^p_{m-2\theta p}}\,dt \leq C\|f\|_{L^p_{m}}\int_0^\infty e^{-t}t^{-\theta}\,dt\, 
 \end{align*}
which yields $D(B_{m,p}^d) \subset D$ where 
$D= \{
u\in D(B_{m,p, max}): \;
 y^{-2\theta}u\in L^p_m \}$
 for  every $0\leq\theta <1$  such  that $ c-1+2\theta<\frac{m+1}{p} <2 $. The equality $D(B_{m,p}^d) = D$ follows from the injectivity of $I-B$ on $D$, as in Proposition \ref{dominioN}, since the function $u_2(y)=y^{\frac{1-c}{2}}K_{\frac{1-c}{2}}\approx c\ \neq 0$ does not belong to $D$ (choosing $2\theta$ sufficiently close to $\frac{m+1}{p}-(c-1)$ or to $2$).

To reach the case when $\theta=1$ and to add the integrability of $D_yu, D_{yy}u$ we argue as in the proposition above. If $g=Bu$ then
$$D_y u(y)=y^{-c}\int_1^y g(s)s^c\, ds +Ky^{-c}:=w(y)+Ky^{-c}.$$ H\"older inequality gives $|w(y)| \le C \|g\|_{L^p_m}(y^{-c}+y^{1-\frac{m+1}{p}})$  and then the assumption $c<1$ and $(m+1)/p<2$ show that $D_yu \in L^1(0,1)$ (with respect to the Lebesgue measure). It follows that $\lim_{y \to 0}u(y)$ exists finite and must be 0, by the same argument for $u_2$. Then $|u(y)| \le C \|g\|_{L^p_m}(y^{1-c}+y^{2-\frac{m+1}{p}})$. At this point the estimate for  $(1\wedge y)^{-2\theta+1}D_y u$ is elementary when $\theta <1$ (and that for $D_{yy }u$ follows from the equation, multiplying by $y^{2-2\theta}$). If $\theta=1$, that is when $c+1<(m+1)/p<2$, then $Ky^{-c-1} \in L^p_m(0,1)$ and $\frac{w}{y} \in L^p_m$ by Hardy inequality, see Lemma \ref{Hardy}. The integrability of $\frac{u}{y^2}$ is proved as in Corollary \ref{dominioN1} (iii). Using $u(0)=0$ we obtain 
$$
\frac{|u(y)|}{y^2} \le  \int_0^1 \eta\, d\eta \int_1^\infty |g(t\eta y)| t^c\, dt+K y^{-c-1}.
$$
Since $y^{-c-1} \in L^p_m(0,1)$ we may assume that $K=0$ 
and Minkowski inequality gives
$$
\left \|\frac{u}{y^2} \right \|_{L^p_m} \le \int_0^1 \eta\, d\eta \int_1^\infty t^c \|g(t\eta \cdot)\|_{L^p_m}\, dt=\|g\|_{L^p_m} \int_0^1 \eta^{1-\frac{m+1}{p}} d \eta \int_1^\infty t^{c-\frac{m+1}{p}} dt.
$$
\qed

Note that when $c<(m+1)/p <2$, $\theta =1/2$ is allowed and $D(B^d_{m,p}) \subset W^{1,p}_m$. The embeddings above do not hold outside the indicated ranges: just take $u(y)=y^{1-c}$ near $0$.

\smallskip
For the next corollary we recall that $\lim_{y \to 0}u(y)$ and $\lim_{y \to 0}D_yu(y)$ exists finite if $u \in W^{2,p}_m$ when $(m+1)/p<1$, see Lemma \ref{l1}, and that both are equal to $0$, when $m \le -1$, see Lemma \ref{mle-1}. 

When $1 \le (m+1)/p <2$, H\"older inequality gives $|D_y u| \le C\|D_{yy} u\|_{L^p_m}(1+y^{1-\frac{m+1}{p}})$ if $1<(m+1)/p <2$,  or  $|D_y u| \le C\|D_{yy} u\|_{L^p_m}(1+|\log y|)$ if $m=p-1$. In both cases $D_y u \in L^1(0,1)$ and $u(0)$ exists finite.

\begin{cor} \label{Dir1} Assume that $c<1$ and $c+1<\frac{m+1}{p}<2$. Then $D(B_{m,p}^d) \subset W^{2,p}_m$ and 
\begin{itemize}
\item[(i)] If $m \le -1$, then $D(B_{m,p}^d)=W^{2,p}_m;$
\item[(ii)]If $0<\frac{m+1}{p}<1$ then $D(B_{m,p}^d)=\{u \in W^{2,p}_m: u(0)=D_yu(0)=0 \}$
\item[(iii)] If $1 \leq \frac{m+1}{p} <2$, then $D(B_{m,p}^d)=\{u \in W^{2,p}_m: u(0)=0\}$.
\end{itemize}
\end{cor}
{\sc Proof. } By Proposition \ref{Dirichlet} the inclusion $D(B_{m,p}^d) \subset \{u \in W^{2,p}_m: u(0)=0\}$ follows. When $(m+1)/p<1$, then also $D_y u(0)=0$ otherwise $D_y u/y \not \in L^P_m$. This shows that in all cases $D(B_{m,p}^d) $ is contained in the right hand sides. To show the equality it suffices, therefore, to note that the function  $u_2(y)=y^{\frac{1-c}{2}}K_{\frac{1-c}{2}}\approx c\ \neq 0$ (see the proof of Proposition \ref{Dirichlet})  does not belong to the right  hand sides (in case (iii) observe also that $D_y u \in L^1(0,1)$ so that $u(0)$ exists finite).
\qed

\section{Degenerate operators in weighted  spaces }\label{Section Degenerate}
In this section we add a potential term to $B$ and study  the whole operator 
$$
L=D_{yy}+\frac{c}{y}D_y-\frac{b}{y^2}
$$ on the (open) half line $\R_+=]0, \infty[$. However, we shall consider $L$ only with Dirichlet boundary conditions at $0$, hence $L=B^d$, when $b=0$, with the understanding that $B^n=B^d$ when $c \ge 1$.


If $1<p<\infty$, we define the maximal operator $L_{p,max}$ through the domain
\begin{equation} \label{Lmax}
D(L_{m,p,max})=\{u \in L_m^p \cap W^{2,p}_{loc} (\R_+ ): Lu \in L_m^p\}.
\end{equation}

The equation $Lu=0$ has solutions $y^{-s_1}$, $y^{-s_2}$ where $s_1,s_2$ are the roots of the indicial equation $f(s)=-s^2+(c-1)s+b=0$ given by

\begin{equation} \label{defs}
s_1:=\frac{c-1}{2}-\sqrt{D},
\quad
s_2:=\frac{c-1}{2}+\sqrt{D}
\end{equation}
where
\bigskip
\begin{equation} \label{defD}
D:=
b+\left(\frac{c-1}{2}\right)^2.
\end{equation}

The above numbers are real if and only if $D \ge 0$. When $D<0$ the equation $u-Lu=f$ cannot have positive distributional solutions for certain positive $f$, see \cite{met-soba-spi3}.
    When $b=0$, then $\sqrt D=|c-1|/2$ and $s_1=0, s_2=c-1$ for $c \ge 1$ and $s_1=c-1, s_2=0$ for $c<1$.

\medskip

%

A multiplication operator transforms $L$  into a  Bessel operator and allows  to transfer the results of the previous section to this more general situation.

\begin{lem}\label{Isometry action} For $k\in\R$, let  
$
(T_k u)(y):=y^ku(y), y>0.
$
Then $T_k$ maps isometrically  $L^p_{m+kp}$ onto $L_{m}^p$ and 
for every $u\in W^{2,1}_{loc}\left(\R_+\right)$ one has
$$
T_{-k} LT_ku= \tilde L u:=D_{yy}u +\frac{\tilde c}{y}D_y u-\frac{\tilde b}{y^2}u
$$

\begin{equation}
\label{tilde b}
\tilde b=b-k\left(c+k-1\right),\qquad \tilde c=c+2k.
\end{equation}

Moreover the discriminant $\tilde D$ and the parameter $\tilde \gamma$, $\tilde s_{1,2}$ of $\tilde L$ defined  in \eqref{defD}, \eqref{defs}  are given by
\begin{align}\label{tilde s gamma}
\tilde D=D,\quad \tilde s_{1,2}=s_{1,2}+k,\quad {\tilde s}^\ast_{1,2}= s^\ast_{1,2}-k.
\end{align}
\end{lem}

Observe now that, choosing $k=-s_i$, $i=1,2$, we get  $\tilde b=0$, $\tilde c_i=c-2s_i$. The operators  $$T_{s_i}LT_{-s_i}=B_i=D_{yy}+\frac{c-2s_i}{y}D_y$$ are therefore  Bessel operators.  The following two results can be found  also in  \cite{met-negro-spina 5} and \cite{met-soba-spi3} when $m=0$ or $m=c$.



\begin{prop}\label{gen-dom-subcritBDirichlet1}
Assume that $D> 0$. If $s_1<\frac{m+1}{p}<s_2+2$  then $L_{m,p}=(L,  D(L_{m,p}))$
where
\begin{align} \label{dd}
\nonumber D(L_{m,p})&=
\bigg\{
u\in D(L_{m,p,max})\;:\;
y^{-2\theta}u\in L_m^p\quad\\& \textrm{for every }\ \theta\in (0,1]\ \textrm{such that}\ s_1+2\theta<\frac{m+1}{p}<s_2+2.
\bigg\}
\end{align}
generates a bounded positive analytic semigroup of angle $\pi/2$ on $L_m^p$. Moreover
\begin{itemize}
\item[(i)] $(1\wedge y)^{2-2\theta}D_{yy}u, (1\wedge y)^{1-2\theta}D_y u \in L^p_m$ for every $u \in D(L_{m,p})$ and $\theta$ as above;
\item[(ii)]
$\lim_{y \to 0}y^{s_2}u(y)=0$ for every $u \in D(L_{m,p}^d)$.
\end{itemize}
 \end{prop}
{\sc Proof.} We use the identity  $T_{s_2}LT_{-s_2}=D_{yy}+\frac{c-2s_2}{y}D_y:=B^d$ and apply Proposition \ref{Dirichlet} in $L^p_{m-s_2p}$. Note that $c-2s_2=1-2\sqrt{D}< 1$ and that 
 $s_1+2 \theta <\frac{m+1}{p}<s_2+2$ is equivalent to $c-2s_2 -1+2\theta <\frac{m-s_2p+1}{p}<2$. Since, by definition, $D(L_{m,p})=T_{-s_2}D(B^d_{m-s_2p,p})$, \eqref{dd} is immediate. The verification of $(i)$ and $(ii)$ is similar. If $u=y^{-s_2} v \in D(L_{m,p})$, then $\lim_{y \to 0} y^{s_2}u(y)=\lim_{y \to 0}v(y)=0$ and $y^{1-2\theta}D_y u=y^{-s_2}(y^{1-2\theta}D_y v-s_2 y^{-2\theta}v) \in L^p_m$, by Proposition \ref{Dirichlet}, again.
\qed

\medskip

Let us now turn to the case $D=0$, where $s_1=s_2$.

\begin{prop}\label{gen-dom-critBDirichlet1}
Assume that $D= 0$. If $s_1<\frac{m+1}{p} <s_1+2$  then $L_{m,p}=(L,  D(L_{m,p}))$
where 
\begin{align*}
 D(L_{m,p})&=
\bigg\{
u\in D(L_{m,p,max})\;:\;  \exists \lim_{y \to 0}y^{s_1}u(y) \in \C
\bigg\}\\[1ex]
&=\left\{
u\in D(L_{m,p,max})\;:\;
y^{-2\theta_0}|\log y|^{-\frac{2}{p}}u
\in L^p_m \Bigl (0, \frac{1}{2} \Bigr)
\right\}
\end{align*}
with $\theta_0=\frac{1}{2}(\frac{m+1}{p}-s_1)\in (0,1)$
generates a bounded positive analytic semigroup of angle $\pi/2$ on $L_m^p$. 
Moreover, 
$$(1\wedge y)^{2-2\theta}D_{yy}u, (1\wedge y)^{1-2\theta}D_y u \in L^p_m$$ 
for every $u \in D(L_{m,p})$ and $s_1+2\theta <\frac{m+1}{p}< s_1+2$.


 \end{prop}
{\sc Proof.} Let us write
\begin{align*}
D_1\!&=\!  \{u\in D(L_{m,p,max}): \lim_{y \to 0}y^{s_1}u(y) \in \C\},
\\[1ex]
 D_2\!&=\! \{u\in D(L_{m,p,max}): y^{-2\theta_0}|\log y|^{-\frac{2}{p}}u
\in L^p_m  (0, \frac{1}{2} )
\}
\end{align*}
and note that $D_1 \subset D_2$, by the choice of $\theta_0$.

We use the identity  $T_{s_1}LT_{-s_1}=D_{yy}+\frac{c-2s_1}{y}D_y=D_{yy}+\frac{1}{y}D_y=B^n$ and apply Proposition \ref{dominioN} in $L^p_{m-s_1p}$ since $c-2s_1=1-2\sqrt{D}= 1$. Note  that 
 $s_1+2 \theta <\frac{m+1}{p}<s_1+2$ is equivalent to $2\theta <\frac{m-s_1p+1}{p}<2$.  If $u=y^{-s_1} v \in D(L_{m,p})=T_{-s_1}D(B^n_{m-s_1p,p})$,  then $\lim_{y \to 0} y^{s_1}u(y)=\lim_{y \to 0}v(y) \in \C$, by Corollary \ref{dominioN1} (i) (and similarly $y^{1-2\theta}D_y u=y^{-s_1}(y^{1-2\theta}D_y v-s_1 y^{-2\theta}v) \in L^p_m$). 
This shows that $D(L_{m,p})\subset D_1 \subset D_2$ and the equality follows since $I-L$ is injective on $D_2$. In fact, if $u-Lu=0$, then $v=y^{s_1}u$ solves $v-Bv=0$, hence $v(y)=c_1 I_0+c_2 K_0$. However, $c_1=0$, since $I_0$ grows exponentially at $\infty$ and $c_2=0$ since $K_0 \approx -\log y$, as $y \to 0$, hence does not satisfy the integrability condition required by $D_2$ near $y=0$.
\qed

An alternative description of the domain is contained in the following proposition, where we do not need to distinguish between $D>0$ and $D=0$.
\begin{prop}\label{gen-dom-subcritBDirichlet1 cor}
If  $D\geq 0$ and  $s_1<\frac{m+1}{p}< s_2+2$,  then 
\begin{align} \label{dd cor}
D(L_{m,p})&=
\bigg\{
u\in D(L_{m,p,max})\;:\;s_1\frac{u}{y^2}+\frac{D_y u}{y}\in L^p_m
\bigg\}.
\end{align}
 \end{prop}
{\sc Proof.} As in the case $D=0$
we use the identity $T_{s_1}LT_{-s_1}=D_{yy}+\frac{c-2s_1}{y}D_y=\tilde B^n$ (this last is in $L^p_{m-s_1p}$), after observing that $c-2s_1=1+2\sqrt{D}\ge 1$. Note  that the conditions  $s_1+2 \theta <\frac{m+1}{p}<s_2+2$ and $2\theta<\frac{m-s_1p+1}{p}<c-2s_1+1$  are equivalent. 

This definition yields the same operator as in Proposition \ref{gen-dom-subcritBDirichlet1} if (and only if)
 $T_{-s_2}D(B^d_{m-s_2p,p})=T_{-s_1}D(\tilde B^n_{m-s_1p,p})$ but this holds since $L$ endowed with both domains generates a semigroup and the first contains the second.

 Indeed, let $u\in T_{s_2-s_1}D(\tilde B^n_{m-s_1p,p})$, that is  $y^{s_1-s_2}u\in  D(\tilde B^n_{m-s_1p,p})$. Then by construction $\tilde B^n_{m-s_1p,p}(y^{s_1-s_2}u)\in L^p_{m-s_1p}$ is equivalent to $B^d_{m-s_2p,p}u\in L^p_{m-s_2p}$. Analogously,  Corollary \ref{dominioN1} applied to $y^{s_1-s_2}u$ yields $y^{-2\theta }u\in L^p_{m-s_2p}$ for every $0\leq\theta \le 1$ such  that $c-s_2-1+2\theta<\frac{m-s_2p+1}{p} <2$. By Proposition \ref{gen-dom-subcritBDirichlet1} this proves that $u\in D(\tilde B^n_{m-s_1p,p})$ i.e. $T_{s_2-s_1}D(\tilde B^n_{m-s_1p,p})\subseteq D(\tilde B^n_{m-s_1p,p})$.

Applying  now  Proposition  \ref{dominioN} and Corollary \ref{dominioN1} to $v=y^{s_1}u$ we get, in addition, that $u\in D(L_{m,p})$ if and only if  $u\in D(L_{m,p,max})$ and 
$$
(a) \quad s_1\frac{u}{y^2}+\frac{D_y u}{y}\in L^p_m, \quad
(b) \quad s_1(s_1-1)\frac{u}{y^2}+2s_1\frac{D_y u}{y}+D_{yy}u\in L^p_m.
$$
However, $(b)$ follows from $(a)$ and $u \in D(L_{m,p,max})$ since  
\begin{align*}
Lu-\left(s_1(s_1-1)\frac{u}{y^2}+2s_1\frac{D_y u}{y}+D_{yy}u\right)&=\left((2-c)s_1-2b\right)\frac{u}{y^2}+(1+2\sqrt D)\frac{D_y u}{y}\\[1ex]
&=(1+2\sqrt D)\left[s_1\frac{u}{y^2}+\frac{D_y u}{y}\right].
\end{align*}
\qed

\medskip

\medskip
We now deduce the estimates for the heat kernel and its derivative for the operator $L$ from those for $B$. We shall consider only the case $s_1 \neq 0$; in fact, if $s_1=0$, then $b=0$ and $c \ge 1$, hence $L=B^n=B^d$ and the estimates are those of Propositions \ref{Estimates Bessel kernels}, \ref{Estimates gradient kernel bessel}.
\begin{prop}\label{kernelL} Let  $1<p <\infty$ such that $0 \neq s_1<\frac{m+1}{p}<s_2+2$.  Then for $z \in \C_+$
\begin{align}\label{Int repr m}
e^{zL}f(y)=\int_{\R^+}p_L(z,y,\rho)f(\rho) d\rho, \quad f\in L_m^p 
\end{align}
where
$$p_L(z,y,\rho)=\frac{1}{2z}y^{\frac{1-c}2}\rho^{\frac{1+c}2}\exp\left\{-\frac{y^2+\rho^2}{4z}\right\}I_{\sqrt D}\left(\frac{y\rho}{2z}\right).$$
For every $\eps>0$, there exist $C_\eps>0$ and $\kappa_{\eps}>0$ such that 
for $z\in\Sigma_{\frac \pi 2-\eps}$
$$
| p_L(z,y,\rho)|\leq C_\eps |z|^{-\frac{1}{2}} \left (\frac{y}{|z|^{\frac{1}{2}}}\wedge 1 \right)^{-s_1} \left (\frac{\rho}{|z|^{\frac{1}{2}}}\wedge 1 \right)^{-s_1+c}
\exp\left(-\frac{|y-\rho|^2}{\kappa_\eps |z|}\right)
$$
and
$$
y^{-1} |p_L(z,y,\rho)|+ |D_{y}p_L(z,y,\rho)|\leq C_\eps |z|^{-1} \left (\frac{y}{|z|^{\frac{1}{2}}}\wedge 1 \right)^{-s_1-1} \left (\frac{\rho}{|z|^{\frac{1}{2}}}\wedge 1 \right)^{-s_1+c}
\exp\left(-\frac{|y-\rho|^2}{\kappa_\eps |z|}\right).
$$

\end{prop}
{\sc{Proof.}}  Let us consider the isometry  $T_{-s_1}:L^p_m\to L^p_{m-s_1p}$. Then for every $u\in L^p_m$
$$e^{zL}u=T_{-s_1}\,e^{zB}\left( T_{s_1}u\right)$$
where $B$ is the pure Bessel  operator $D_{yy}+\frac{c-2s_1}{y}$ on $L^p_{m-s_1p}$.
Let $p_L(z,y,\rho)$ be the heat kernel of $\left(e^{zL}\right)_{z\in\C_+}$ on $L_m^p$. The last relation between the semigroups translate into the analogous equality for the heat kernels:
\begin{align*}
 p_L(z,y,\rho)=y^{-s_1}p_{B^n}(z,y,\rho)\rho^{s_1}.
 \end{align*}
From Proposition \ref{Estimates Bessel kernels} and from Lemma \ref{equiv1} it follows that
\begin{align*}
| p_L(z,y,\rho)|&=y^{-s_1}|p_B^n(z,y,\rho)|\rho^{s_1}\leq \frac{C_\eps}{|z|^\frac{1}{2}} y^{-s_1} \rho^{s_1} \left (\frac{\rho}{|z|^{\frac{1}{2}}}\wedge 1 \right)^{c-2s_1}
\exp\left(-\frac{|y-\rho|^2}{\kappa_\eps |z| }\right)\\&=\frac{C_\eps}{|z|^\frac{1}{2}}  \left(\frac{y}{|z|^\frac{1}{2}}\right)^{-s_1} \left(\frac{\rho}{|z|^\frac{1}{2}}\right)^{s_1} \left (\frac{\rho}{|z|^{\frac{1}{2}}}\wedge 1 \right)^{c-2s_1}
\exp\left(-\frac{|y-\rho|^2}{\kappa_\eps |z| }\right)
\\&\leq \frac{C_\eps}{|z|^\frac{1}{2}}  \left(\frac{y}{|z|^{\frac{1}{2}}}\wedge 1\right)^{-s_1}  \left (\frac{\rho}{|z|^{\frac{1}{2}}}\wedge 1 \right)^{c-s_1}
\exp\left(-\frac{|y-\rho|^2}{\kappa_\eps |z| }\right).
 \end{align*}
Concerning the derivative with respect to $y$, we get 
\begin{align*}
 D_yp_L(z,y,\rho)=-s_1y^{-s_1-1}p_{B^n}(z,y,\rho)\rho^{s_1}+y^{-s_1}D_yp_{B^n}(z,y,\rho)\rho^{s_1}.
 \end{align*}
From Proposition \ref{Estimates gradient kernel bessel}, Proposition \ref{Estimates Bessel kernels} and  Lemma \ref{equiv1} it follows that
\begin{align*}
 D_yp_L(z,y,\rho)&=-s_1y^{-s_1-1}p_{B^n}(z,y,\rho)\rho^{s_1}+y^{-s_1}D_yp_{B^n}(z,y,\rho)\rho^{s_1}\\&\leq 
\frac{C_\eps}{|z|}  \left(\frac{y}{|z|^\frac{1}{2}}\wedge 1\right)^{-s_1-1} \left (\frac{\rho}{|z|^{\frac{1}{2}}}\wedge 1 \right)^{c-s_1}
\exp\left(-\frac{|y-\rho|^2}{\kappa_\eps |z| }\right).
 \end{align*}
The estimate for $y^{-1}p_L$ follows easily from that of $p_L$.
\qed

\section{Remarks on domain characterization and uniqueness}
The domain characterizations for $B^n, B^d, L$ can be stated  by adding explicit estimates. For example,  in Proposition \ref{dominioN}, one can add 
$$\|D_{yy}u\|_{L^p_m} +\| y^{-1}D_y u\|_{L^p_m} \le C\|Bu\|_{L^p_m}, \quad u \in D(B^n_{m,p})$$ (the additional term $\|u\|_{L^p_m}$ does not appear, by scaling). This follows from the proof (actually, this is the proof) but can be also deduced from the statement by the closed graph theorem. This remark applies to all the domain characterizations including those of the next sections; we decided not to write down them explicitly for exposition reasons  but we shall use in Section 8. 

\smallskip

As already pointed out, the assumption $D \ge 0$ is crucial for positivity and always satisfied in the case of Bessel operators. When $D<0$ the equation $u-Lu=f$ cannot have positive distributional solutions for certain positive $f$, see \cite{met-soba-spi3}. However, $L$ can be the generator of a semigroup even when $D<0$, see \cite{met-soba} for Schr\"odinger operators with inverse square potential $\Delta-b|y|^{-2}$ with $b<-1/4$.
\\

We note that the results for $B^d$ can be deduced from those for $\tilde B^n$ (here we denote by $B, \tilde B$ two different but related  Bessel operators). This simple (but surprising) fact is actually the core of the proof of  Proposition \ref{gen-dom-subcritBDirichlet1 cor} where the operators $L_{m,p}$ are transformed via similarity to  pure Bessel operators with  $c \geq 1$. This approach, however, needs a change in the reference measures and the knowledge of Bessel operators in $L^p_m$ for every admissible $m$. We prefer to start with a form in $L^2$ of the symmetrizing measure, as it is usually done. Moreover, using both the direct approach and the transformation above, one gets different and complementing descriptions of $D(L_{m,p})$, see Propositions \ref{gen-dom-subcritBDirichlet1 cor}, \ref{gen-dom-subcritBDirichlet1} and subsequent corollaries,  which are likely difficult to discover simultaneously using only one method. \\

A natural question arises  if different boundary conditions can be imposed to produce different semigroups. This is the case, for example, for the Bessel operators of Section 3, in the range $-1<c<1$ where $B^n \neq B^d$.
To state the uniqueness question more precisely we define
$L_{m,p,min}$ as the closure, in $L^p_m$ of $(L,C_c^\infty (\R_+))$ (the closure exists since this operator is contained in the closed operator $L_{m,p,max}$) and it is clear that $L_{m,p,min} \subset L_{m,p,max}$. We look at realizations  $L_D=(L,D)$, such that $L_{m,p,min} \subset L_D \subset L_{m,p,max}$. The following results can be found in \cite[Propositions 2.4, 2.5]{met-negro-spina 7} and \cite[Propositions 3.12, 3.28, 3.30]{met-soba-spi3} in the $N$-dimensional case and for $m=0$. The generalization to any $m$ is straightforward, through the transformation $T_k$.
\begin{prop} \label{domainL}
\begin{itemize}
\item[(i)] If $\frac{m+1}{p} \not \in (s_1, s_2+2)$, then no realization  $L_{m,p,min} \subset L_D \subset L_{m,p,max}$ generates a semigroup in $L^p_m$;
\item[(ii)] $L_{m,p,max}$ generates a semigroup if and only if $s_1<\frac{m+1}{p} \leq  s_2$;
\item[(iii)]  $L_{m,p,min}$ generates a semigroup if and only if $s_1+2 \leq \frac{m+1}{p}< s_2+2$.
\end{itemize}
\end{prop}
%
%

In particular $L$ generates a unique semigroup in cases $(ii)$ and $(iii)$ and $L_{m,p}= L_{m,p,max}$ or $L_{m,p}= L_{m,p,min}$, respectively. 

Therefore if the intervals $(s_1,s_2]$ and $[s_1+2,s_2+2)$ overlap, that is if 
$s_1+2 \leq s_2$ or equivalently $D\geq 1$,
we have uniqueness in all $L^p_m$ for which there is generation and, moreover, $L_{m,p,max}=L_{m,p,min}$ if $(m+1)/p \in [s_1+2,s_2]$. 

Uniqueness fails if
$ s_2<s_1+2 $, i.e. $0\leq D<1$, and $ (m+1)/p \in (s_2,s_1+2)$, as we show below but only for $D>0$.

\begin{prop}
If  $0< D<1$  and $s_2<\frac{m+1}p<s_1+2$,   then $(L,  D(L))$
where
\begin{align} \label{L neumann}
D(L)&=
\bigg\{
u\in D(L_{m,p,max})\;:\;s_2\frac{u}{y^2}+\frac{D_y u}{y}\in L^p_m
\bigg\}
\end{align}
generates  a bounded  positive analytic semigroup of angle $\pi/2$ on $L_m^p$.
 \end{prop}
{\sc Proof.} We proceed as  in the proof of  Proposition \ref{gen-dom-subcritBDirichlet1 cor} but in place of the isometry $T_{-s_1}$  we use  the identity $T_{s_2}LT_{-s_2}=D_{yy}+\frac{c-2s_2}{y}D_y=\tilde B^n$ (this last is in $L^p_{m-s_2p}$), after observing that, under the given hypotheses,   $c-2s_2=1-2\sqrt{D}>-1$. Note  that the conditions  $s_2  <\frac{m+1}{p}<s_1+2$ and $0<\frac{m-s_2p+1}{p}<c-2s_2+1$  are equivalent.  The generation result follows  then by similarity from  Proposition \ref{Gen DN}. The description of $ D(L)$ follows by applying  Proposition \ref{dominioN} to $v=y^{s_2}u\in D(\tilde B^n_{m-s_2p,p})$, $u \in D(L)$,   as  in the proof of  Proposition \ref{gen-dom-subcritBDirichlet1 cor}.\qed

 We point out that in the range $s_2<\frac{m+1}p<s_1+2$ the operators $L_{m,p}$ of Section 4 and $(L,D(L))$ just constructed are different. In fact, a compactly supported function  $u$ which is equal to $y^{-s_1}$ in a neighborhood of $0$, belongs to $D(L_{m,p})$ but not to $D(L)$ (otherwise $y^{-2}u$ would be in $L^p_m$ and this is not the case since $(m+1)/p <s_1+2$).

 When $c=0$, then $L$ is a 1d-Schr\"odinger operator with inverse square potential  and the condition $s_2<s_1+2$ becomes $-\frac 14 \leq b <\frac 34$. That uniqueness really does not occur is proved for example in \cite{met-soba-spi2} where different positive and analytic semigroups are exhibited in $L^2$. In the case of Bessel operators, the condition $s_2<s_1+2$ becomes $-1<c<3$, and not $-1<c<1$ as one could guess.  

We close this section by describing cores for these degenerate operators. Observe that by (iii) of the above proposition, $C_c^\infty (0,\infty)$ is a core for $L_{m,p}$ if and only if $ s_1+2 \leq (m+1)/p <s_2+2$. 
\begin{prop} \label{core}
\begin{itemize}
\item[(i)]
If $c>-1$ and $0<\frac{m+1}{p}<c+1$, then $\mathcal {D}= \left\{u \in C_c^\infty ([0,\infty)): D_y u(0)=0  \right \}$ is a core for $B^n_{m,p}$.
\item[(ii)]
 If $s_1<\frac{m+1}{p}<s_2+2$, then $\mathcal D=\left \{u=y^{-s_1}v: v \in C_c^\infty ([0, \infty), \ D_y v(0)=0 \right \}$ is a core for $D(L_{m,p})$.
\end{itemize}
\end{prop}
{\sc Proof. } (i) By Proposition \ref{dominioN}, $\mathcal D \subset D(B^n_{m,p})$. Let $u \in D(B^n_{m,p})$, $f=(I-B^n)u $, $f_k=f \chi_{[k^{-1}, \infty)}$ and $u_k=(I-B^n_{m,p})^{-1}f_k$ so that  $u_k \to u$ with respect to the graph norm. By Proposition \ref{risolventeB^n}, $u_k(y)=c_k y^{\frac{1-c}{2}}I_{\frac{c-1}{2}}(y)$ for $y \leq \frac 1k$, hence $u_k \in C^\infty ([0,\frac 1k]) $ and  $D_y u_k(0)=0$ by Lemma \ref{behave}. Now the proof is straightforward. We fix $k$ and a cut-off function  $\phi $ which is equal to $1$ in $[0 \frac{1}{2k}]$ and to $0$ for $y \geq \frac 1k$; then  write $u_k=\phi u_k +(1-\phi)u_k$ and smooth $(1-\phi)u_k$ by using convolutions (plus cut-off at $\infty$ to make everything with compact support).

The proof of (ii) is similar using the Green function of $L_{m,p}$ or  using the transormation $T_{-s_1}$ as in Proposition \ref{gen-dom-subcritBDirichlet1 cor}.
\qed

\section{Sums of closed  operators}
The operator $\mathcal L=\Delta_x+L_y$ (we write $\Delta_x$, $L_y$ to indicate the variables on which the operators act) is the sum of the Laplacian $\Delta_x$ and of the degenerate  one dimensional operator $L_y$ which clearly commute on smooth functions. Regularity properties for $\mathcal L$ follow once we prove the estimate
\begin{equation}  \label{regularity}
\|\Delta_x u\|_p+\|L_y u\|_p \le C \| \mathcal L u\|_p
\end{equation}
 where the $L^p$ norms are taken over $\R^{N+1}_+$ on a sufficiently large set of functions $u$. This is equivalent to saying that the domain of $\mathcal L $ is the intersection of the domain of $\Delta_x$ and $L_y$ (after aprropriate tensorization) or that the operator $\Delta_x {\mathcal L}^{-1}$ is bounded. This strategy arose first in the study of maximal regularity of parabolic problems, that is for the equation $u_t=Au+f, u(0)=0$ where $A$ is the generator of an analytic semigroup on a Banach space $X$. Estimates like
$$
\|u_t\|_p+\|Au\|_p \le \|f\|_p
$$
where now the $L^p$ norm is that of $L^p([0, T[;X)$ can be interpreted as closedness of $D_t-A$ on the intersection of the respective domains or, equivalently, boundedness of the operator $A(D_t-A)^{-1}$ in $L^p([0, T[;X)$.

Nowadays this strategy is well established and relies on Mikhlin vector-valued multiplier theorems.
Let us state the relevant definitions and main results we need, referring the reader to \cite{DHP}, \cite{PS} or \cite{KW}.

Let ${\cal S}$ be a subset of $B(X)$, the space of all bounded linear operators on a Banach space $X$. ${\cal S}$ is $\mathcal R$-bounded if there is a constant $C$ such that
$$
\|\sum_i \eps_i S_i x_i\|_{L^p(\Omega; X)} \le C\|\sum_i \eps_i  x_i\|_{L^p(\Omega; X)} 
$$
for every finite sum as above, where $(x_i ) \subset X, (S_i) \subset {\cal S}$ and $\eps_i:\Omega \to \{-1,1\}$ are independent and symmetric random variables on a probability space $\Omega$.
It is well-known that this definition is independent of $1 \le p<\infty$ and that $\mathcal R$-boundedness is equivalent to boundedness when $X$ is an Hilbert space.
When $X$ is an $L^p$ space (with respect to any $\sigma$-finite measure), testing  $\mathcal R$-boundedness is equivalent to proving square function estimates, see \cite[Remark 2.9 ]{KW}.

\begin{prop}\label{Square funct R-bound} Let ${\cal S} \subset B(L^p(\Sigma))$, $1<p<\infty$. Then ${\cal S}$ is $\mathcal R$-bounded if and only if there is a constant $C>0$ such that for every finite family $(f_i)\in L^p(\Sigma), (S_i) \in {\cal S}$
$$
\left\|\left (\sum_i |S_if_i|^2\right )^{\frac{1}{2}}\right\|_{L^p(\Sigma)} \le C\left\|\left (\sum_i |f_i|^2\right)^{\frac{1}{2}}\right\|_{L^p(\Sigma)}.
$$
\end{prop}
By the proposition above $\mathcal R$-boundedness follows from domination. We formualate this simple but important fact as  a corollary.
\begin{cor} \label{domination}
Let  ${\cal S}, {\cal T} \subset B(L^p(\Sigma))$, $1<p<\infty$ and asuume that $\cal T$ is $\mathcal R$ bounded and that for every $S \in \cal S$ there exists $T \in \cal T$ such that $|Sf| \leq |Tf|$ pointwise, for every $f \in L^p(\Sigma)$. Then ${\cal S}$ is $\mathcal R$-bounded.
\end{cor}
{\sc Proof. } This follows since 
$$\left (\sum_i |S_if_i|^2\right )^{\frac{1}{2}} \leq \left (\sum_i |T_if_i|^2\right )^{\frac{1}{2}}
$$
pointwise.
\qed

Let $(A, D(A))$ be a sectorial operator in a Banach space $X$; this means that $\rho (-A) \supset \Sigma_{\pi-\phi}$ for some $\phi <\pi$ and that $\lambda (\lambda+A)^{-1}$ is bounded in $\Sigma_{\pi-\phi}$. The infimum of all such $\phi$ is called the spectral angle of $A$ and denoted by $\phi_A$. Note that $-A$ generates an analytic semigroup if and only if $\phi_A<\pi/2$. The definition of $\mathcal  R$-sectorial operator is similar, substituting boundedness of $\lambda(\lambda+A)^{-1}$ with $\mathcal R$-boundedness in $\Sigma_{\pi-\phi}$. As above one denotes by $\phi^R_A$ the infimum of all $\phi$ for which this happens; since $\mathcal R$-boundedness implies boundedness, we have $\phi_A \le \phi^R_A$.

\medskip

The $\mathcal R$-boundedness of the resolvent characterizes the regularity of the associated inhomogeneous parabolic problem, as we explain now.

An analytic semigroup $(e^{-tA})_{t \ge0}$ on a Banach space $X$ with generator $-A$ has
{\it maximal regularity of type $L^q$} ($1<q<\infty$)
if for each $f\in L^q([0,T];X)$ the function
$t\mapsto u(t)=\int_0^te^{-(t-s)A})f(s)\,ds$ belongs to
$W^{1,q}([0,T];X)\cap L^q([0,T];D(B))$.
This means that the mild solution of the evolution equation
$$u'(t)+Au(t)=f(t), \quad t>0, \qquad u(0)=0,$$
is in fact a strong solution and has the best regularity one can expect.
It is known that this property does not depend on $1<q<\infty$ and $T>0$.
A characterization of maximal regularity is available in UMD Banach spaces, through the $\mathcal  R$-boundedness of the resolvent in sector larger than the right half plane or, equivalently, of the semigroup in a sector around the positive axis. In the case of $L^p$ spaces it can be restated in the following form,  see \cite[Theorem 1.11]{KW}

\begin{teo}\label{MR} Let $(e^{-tA})_{t \ge0}$ be a bounded analytic semigroup in $L^p(\Sigma)$ with generator $-A$. Then $T(\cdot)$ has maximal regularity of type $L^q$  if and only if there are constants $0<\phi<\pi/2 $, $C>0$ such that for every finite sequence $(\lambda_i) \subset \Sigma_{\pi/2+\phi}$, $(f_i) \subset  L^p$
$$
\left\|\left (\sum_i |\lambda_i (\lambda_i+A)^{-1}f_i|^2\right )^{\frac{1}{2}}\right\|_{L^p(\Sigma)} \le C\left\|\left (\sum_i |f_i|^2\right)^{\frac{1}{2}}\right\|_{L^p(\Sigma)}.
$$
or, equivalently for every finite sequence
$(z_i) \subset \Sigma_\phi$, $(f_i) \subset  L^p$
$$
\left\|\left (\sum_i |e^{-z_i A}f_i|^2\right )^{\frac{1}{2}}\right\|_{L^p(\Sigma)} \le C\left\|\left (\sum_i |f_i|^2\right)^{\frac{1}{2}}\right\|_{L^p(\Sigma)}.
$$
\end{teo}
\medskip

Finally we state the operator-valued Mikhlin Fourier multiplier theorem in the N-dimensional case, see \cite[Theorem 3.25]{DHP} or \cite[Theorem 4.6]{KW}.
\begin{teo}   \label{mikhlin}
Let $1<p<\infty$, $M\in C^N(\R^N\setminus \{0\}; B(L^p(\Sigma))$ be such that  the set
$$\left \{|\xi|^{|\alpha|}D^\alpha_\xi M(\xi): \xi\in \R^{N}\setminus\{0\}, \ |\alpha | \leq N \right \}$$
is $\mathcal{R}$-bounded.
Then the operator $T_M={\cal F}^{-1}M {\cal F}$ is bounded in $L^p(\R^N, L^p(\Sigma))$, where $\cal{F}$ denotes the Fourier transform.
\end{teo}

\subsection{Muckenhoupt weighted estimates}

Let  $\left(S,d,\nu\right)$ be a space of homogeneous type, that is a metric space endowed with a Borel measure $\nu$ which is doubling on balls.  When $X=L^p\left(S,d,\nu\right)$  the  square function estimate  in Theorem \ref{Square funct R-bound} can be reduced to a family of Muckenhoupt weighted estimates of the type
$$\|e^{-zA}f\|_{L^p(w)}\leq C \|f\|_{L^p(w)}, \quad z \in \Sigma_\phi, $$
see Theorem \ref{Extrap Mack} below. With this in mind, we recall preliminarily, the definition and the  essential properties about Muckenhoupt weights. For the proof of the following results as well as for further details, we refer the reader to \cite[Chapter 2 and 5]{auscher1} and \cite[Chapter 1]{Stromberg}. Let $w$ be a weight i.e.   a non-negative locally integrable function defined on $S$; we
use the notation
$$
\aver{E} w = \frac {1}{\nu(E)} \int_{E} w(x)\, d\nu(x), \qquad w(E)=\int_{E} w(x)\, d\nu(x).
$$
Let $\mathcal M_\nu$  denote the
uncentered maximal operator over balls  in $S$  defined by 
\begin{align}\label{Maximal operator}
\mathcal M_\nu f(x):=\sup_{B\ni x}\, \aver{B}|f|,\quad x\in S,
\end{align}
where the supremum is taken over all balls  of $S$ containing $x$. We recall that $\mathcal M_\nu$ is bounded on $L^p(w)$ if and only if $w\in A_p$, see for example \cite[Theorem 7.3]{Duo}.

 We say
that $w\in A_p$, $1< p <\infty$, if there exists a constant $C$ such
that for every ball $B\subseteq S$   one has 
\begin{align}\label{Def A_p}
\Big(\aver{B} w\Big)\,
\Big(\aver{B} w^{1-p'}\Big)^{p-1}\le C.
\end{align}
For $p=1$, we say that $w\in A_1$ if there is a constant $C$ such
hat $\mathcal M_\nu w\le C\,w$  a.e..

The weight  $w$ is in the reverse H\"older class of order $q$,  $w\in
RH_{q}$, $1< q \le \infty$, if there is a constant $C$ such that for
every ball $B\subseteq S$
$$
\Big(\aver{B} w^q\Big)^{\frac1q}
\le C\, \aver{B} w,
$$
with the usual modification for $q=\infty$. For $p=1$, $RH_1$ is the set of all weights.
The best constants appearing in the previous inequalities are referred  respectively as the $A_p$ and the $RH_q$ constants of $w$.

We sum up in the following proposition the properties we need about these classes of weights.

\begin{prop}\label{Prop Macke}
The following properties hold:
\begin{itemize}
\item [(i)] $A_1\subset A_p \subset A_q$ for every $1\leq  p\leq q\leq\infty$;
\item [(ii)] $w\in A_p$, $1<p<\infty$, if and only if $w^{1-p'}\in A_{p'}$;
\item [(iii)] If $w\in A_p$, $1<p<\infty$, then there exists $1<q<p$ such that $w\in A_q$;
\item [(iv)] $RH_\infty\subset RH_q\subset RH_p$ for $1 <p\leq q\leq\infty$;

\item [(v)] If $w\in RH_q$, $1<q<\infty$, then there exists $q<p<\infty$ such that $w\in RH_p$;
\item[(vi)] $A_\infty:=\bigcup_{1\le p<\infty} A_p=\bigcup_{1<q\le
\infty} RH_q . $
\item [(vii)] Let $1<p_0<p<q_0<\infty$. Then we have
$$w\in  A_{\frac{p}{p_0}}\cap RH_{\left(\frac{q_0}{p}\right)'}\iff w^{-\frac{p'}{p}}=w^{1-p'}\in A_{\frac{p'}{q_0'}}\cap RH_{\left(\frac{p'_0}{p'}\right)'}.$$
\item[(viii)] If $1\le p\le \infty$ and $1\le r<\infty$ then 
$$
w\in A_p \cap RH_r \quad \Longleftrightarrow \quad w^{r},w^{-\frac 1{p-1}}\in A_{\infty} \quad\Longleftrightarrow\quad w^{r}\in A_{r\,(p-1)+1}.$$   
\end{itemize}
\end{prop}
{\sc{Proof.}} Properties $(i)$-$(vi)$ can be found in \cite[Chapter 7]{Duo}, \cite[Chapter 1]{Stromberg}. Point (vii) follows as in \cite[Lemma 4.4]{auscher1}. The first equivalence in (viii)  is proved in \cite[Lemma 11, Chapter 1]{Stromberg}; the second follows as in \cite{JN}.\qed

A proof of the following result is in  \cite[Corollary 14]{Stromberg} or  \cite[Chapter 7]{Duo}.

\begin{lem}\label{Doubling weight}
Let $w\in A_p\cap RH_r$, $1 < r,p<\infty$. Then there exists a constant $C > 1$ such
that for any ball $B$ and any measurable subset $E\subset B$, 
$$C^{-1}\left(\frac{\nu(E)}{\nu(B)}\right)^p\leq \frac{w(E)}{w(B)}\leq C\left(\frac{\nu(E)}{\nu(B)}\right)^{\frac{r-1}{r}}.$$
\end{lem}

We now state  an  extrapolation result originally due to Rubio de Francia, adapted as in  \cite[Theorem 4.9]{auscher1}, which allows to reduce the  square function estimate  in Theorem \ref{Square funct R-bound}  to a family of Muckenhoupt weighted estimates. 
Only weights and pairs of functions appear and no operator is involved. In what follows we consider  families $\mathcal F=\{(f,g): f, g \in  L_+^0(S) \}$, where $L_+^0(S)$ is the set    
of all  non-negative,
measurable functions  defined on $S$.
\begin{teo}\label{Extrap Mack}
Let  $\left(S,d,\nu\right)$ be a space of homogeneous type and let $\mathcal F \subseteq L_+^0(S)\times L_+^0(S)$.  Suppose that there exists $p$ with $p_0\le
p\le q_0$ (and $p<\infty$ if $q_0=\infty$), such that for $(f,g)\in
\mathcal F$,
\begin{equation*}\label{extrapol-p}
\|f\|_{L^p(w)}
\le
C\|g\|_{L^p(w)},
\qquad
\mbox{for all }w\in A_{\frac{p}{p_0}}\cap
RH_{\left(\frac{q_0}{p}\right)'},
\end{equation*}
 Then, for all $p_0<q<q_0$ and $(f,g)\in\mathcal F$ we have
\begin{equation*}\label{extrapol-q}
\|f\|_{L^q(w)} \le C\,\|g\|_{L^q(w)},
\qquad
\mbox{for all }w\in A_{\frac{q}{p_0}}\cap
RH_{\left(\frac{q_0}{q}\right)'},
\end{equation*}
Moreover, for all $p_0<q,r<q_0$ and $\{(f_j,g_j)\}\subset\mathcal F$ we have
\begin{equation*}\label{extrapol:v-v}
\Big\|
\Big(\sum_j (f_j)^r\Big)^{1/r}
\Big\|_{L^q(w)}
\le
C\,\Big\|
\Big(\sum_j (g_j)^r\Big)^{1/r}
\Big\|_{L^q(w)},
\quad
\mbox{for all }w\in A_{\frac{q}{p_0}}\cap
RH_{\left(\frac{q_0}{q}\right)'}.
\end{equation*}
All the constants  $C$ above may vary from line to line but depend only on the $A_s$ and $RH_s$ constants of $w$.
\end{teo}

Combining Theorem \ref{Square funct R-bound} and Theorem \ref{Extrap Mack} we derive the following characterization of maximal regularity in terms of boundedness over $L^p(w)$ spaces. 
\begin{teo}\label{Maximal regularity Mack}
Let  $\left(S,d,\nu\right)$ be a space of homogeneous type,  $p_0\le
p\le q_0$ with $p<\infty$ if $q_0=\infty$ and  $p_0<2<q_0$.  Let $(e^{-zA})_{z \in \Sigma_\delta}$ be a bounded analytic semigroup in $L^p\left(S,\nu\right)$ defined in a sector $\Sigma_\delta$, $\delta >0$.   Suppose that such that for $f\in L^p\left(S,\nu\right)$,
\begin{equation*}
\|e^{-zA}f\|_{L^p(w)}
\le
C\|f\|_{L^p(w)},\qquad
\mbox{for all }z \in \Sigma_\delta,
\qquad
\mbox{for all }w\in A_{\frac{p}{p_0}}\cap
RH_{\left(\frac{q_0}{p}\right)'},
\end{equation*}
where  $C$ depends only on the $A_s$ and $RH_s$ constants of $w$.

 Then, for all $p_0<q<q_0$,  $(e^{-tA})_{t \ge0}$ has maximal regularity on $L^q(w)$ for all  $w\in A_{\frac{q}{p_0}}\cap
RH_{\left(\frac{q_0}{q}\right)'}$. 
\end{teo}

The following three lemmas will be crucial in the proof of maximal regularity.

\begin{lem}\label{Max functs estimate}
Let   $w\in A_p$, $p \ge 1$,  and  let $\nu_w$ be the measure $wd\nu$. Denote by  $\mathcal{M}_{\nu_w}$ and $\mathcal{M}_{\nu}$ the maximal function defined   by $\nu_w$ and $\nu$. Then $\left(S,d,\nu_w\right)$ is  a space of homogeneous type and 
\begin{align*}
\mathcal{M}_\nu f\leq A_p (w)^\frac 1 p\left(\mathcal{M}_{\nu_w}|f|^p\right)^{\frac 1 p}, \quad f\in L^1_{loc}\left(S,\nu\right),
\end{align*}
where $A_p(w)$ is the $A_p$ constant of $w$.
\end{lem}
{\sc{Proof.}} The doubling condition for the measure $\nu_w$ follows  from that of  $\nu$ and  Lemma \ref{Doubling weight}. To prove the second claim, let $f\in L^1_{loc}\left(S,\nu\right)$. Then for every ball $B$ of $S$ one has, applying H\"older's inequality,   
\begin{align*}
\frac{1}{\nu(B)}\int_B |f|d\nu&=\frac{1}{\nu(B)}\int_B |f|w^{\frac 1 p}w^{-\frac 1 p}d\nu 
\leq \frac{1}{\nu(B)}\left(\int_B |f|^pwd\nu\right)^{\frac 1 p}\left(\int_B w^{1-p'}d\nu\right)^{\frac 1 {p'}}.
\end{align*}
Using \eqref{Def A_p} we get
\begin{align*}
\frac{1}{\nu(B)}\int_B |f|d\nu\leq A_p(w)^\frac 1 p\left(\frac{1}{\nu_w(B)}\int_B |f|^pwd\nu\right)^{\frac 1 p}
\end{align*}
which, taking the supremum over $B$, yields the required claim. The case $p=1$ follows similarly.\qed

\begin{lem}\label{convolution radial}
Let $p$ be a  non-negative, locally integrable function on $\R^M$ and consider the  measure $\nu=p\, dy$. Let $\mathcal{M_\nu}$ be  the uncentered maximal operator relative to $\nu$, defined as  in \eqref{Maximal operator}.  If  $0 \leq \phi\in L^1\left(\R^M,\nu\right)$ is  radial and  decreasing then
\begin{align*}
|(\phi\ast pf)(y)|\leq \|\phi\|_{L^1\left(\R^M,\nu\right)}\mathcal{M}_\nu f(y),\quad y\in\R^M,\quad f\in L^1_{loc}\left(\R^M,\nu\right).
\end{align*}
If $p$ is homogeneous of degree $k$ i.e.  $p(ty)=t^kp(y)$  for all $x\in \R^N$ and $t>0$, then setting $\phi_t:=t^{-M-k}\phi(t^{-1}y)$ one has 
 \begin{align*}
\sup_{ t>0}|(\phi_t\ast pf)(y)|\leq \|\phi\|_{L^1\left(\R^M,\nu\right)}\mathcal{M}_\nu f(y).
\end{align*}
\end{lem}
{\sc{Proof.}}
Let us suppose preliminarily that $\phi$ is a simple function and let us write, for some $a_1,\dots,a_k>0$ and  balls $B_1,\dots, B_k$ centered at $0$,
\begin{align*}
\phi(y)=\sum_{j=1}^k a_j \chi_{B_{j}}(y).
\end{align*}
Then, since $\|\phi\|_{L^1\left(\R^M,\nu\right)}=\sum_{j=1}^k a_j\, \nu(B_{j})$ and $(\chi_{B_{j}}\ast pf)(y)=\int_{y-B_j}f(z)d\nu$ , we get
\begin{align*}
(\phi\ast pf)(y)=\sum_{j=1}^k a_j\, \nu(B_{j})\frac{1}{\nu(B_{j})}(\chi_{B_{j}}\ast pf)(y)\leq \|\phi\|_{L^1\left(\R^M,\nu\right)}\mathcal{M}_\nu f(y).
\end{align*}
In the general case the first required claim follows since $\phi$ can be approximated by a sequence of simple functions which increase to it monotonically. To prove the second claim it is enough to observe that, under the homogeneity assumptions on $p$, one has $\|\phi_t\|_{L^1\left(\R^M,\nu\right)}=\|\phi\|_{L^1\left(\R^M,\nu\right)}$
\qed

\begin{lem}\label{radial weights}
Let $m\in\R$ be such that $M+m>0$ and let $d\mu_m=|y|^m dy$. For every $k\in\R$ let us consider  the radial weight $w(y)=|y|^k$. The following properties hold.
\begin{itemize}
\item[(i)] If $1\leq p\leq\infty$ then  $w\in A_p\left(\mu_m\right)$ if and only if $-(M+m)<k<(M+m)(p-1)$.
\item[(ii)] If $1\le p\le \infty$ and $1\le r<\infty$ then $w\in A_p(\mu_m) \cap RH_r(\mu_m)$ if and only if $-\frac{M+m}{r}<k<(M+m)(p-1)$.
\end{itemize} 
\end{lem}
{\sc{Proof.}} To prove (i), we start by considering balls of center $y_0$ and radius $1$. Fix $R>1$ and assume first that  $|y_0|\leq R$. Then both  $|y|^k$ and $|y|^{-\frac{k}{p-1}}$ are  integrable in $B(y_0,1)$ with respect to the measure $\mu_m$ and 
\begin{equation} \label{Cond-palla}
\left(\frac{1}{\mu_m(B(y_0,1))}\int_{B(y_0,1)}|y|^k\ d\mu_m\right)\left(\frac{1}{\mu_m(B(y_0,1))}\int_{B(y_0,1)}|y|^{-\frac{k}{p-1}}\ d\mu_m\right)^{p-1}\leq C\end{equation}
for some positive constant $C$ depending on $R$. On the other hand, when $|y_0|>R$, then 
$$
\left(\frac{1}{\mu_m(B(y_0,1))}\int_{B(y_0,1)}|y|^k\ d\mu_m\right)\approx |y_0|^k, \quad \left(\frac{1}{\mu_m(B(y_0,1))}\int_{B(y_0,1)}|y|^{-\frac{k}{p-1}}\ d\mu_m\right)^{p-1}\approx |y_0|^{-k}
$$
and the left hand side in \eqref{Cond-palla}
is bounded from above and below by a constant.
For a general ball of radius $r$ the claim follows  by  scaling. Property (ii)   follows using (i) and property (viii) of Proposition \ref{Prop Macke}.  
\qed

\section{$\mathcal{R}$-boundedness for a family of integral operators}

In this section we study the ${\mathcal R}$-boundedness of the  family of integral operators 
\begin{align}\label{gen eq semigrop est Lp_m 1}
S^{\alpha,\beta}(t)f(y)=t^{-\frac {M} 2}\,\left (\frac{|y|}{\sqrt t}\wedge 1 \right)^{-\alpha}\int_{\R^M}  \left (\frac{|z|}{\sqrt t}\wedge 1 \right)^{-\beta}
\exp\left(-\frac{|y-z|^2}{\kappa t}\right)f(z) \,dz,
\end{align}
where $\kappa$ is a positive constant. We omit the dependence on $\kappa$ even though in some proofs we need to vary it.

For $m\in\R$ we consider  the measure $d\mu_m=|y|^m dy$ on $\R^M$ 
and study the action of $S^{\alpha,\beta}(t)$ over  the space $L^p_m=L^p(\R^M, d\mu_m)$ for $1<p<\infty$. We prove   that when $S^{\alpha,\beta}(1)$ is bounded  in $L^p_m$, then the family 
$\left(S^{\alpha,\beta}(t)\right)_{t> 0}$ is  also $\mathcal{R}$-bounded on  $L^p_m$. Note that we write the kernel of these operators always with respect to the Lebesgue measure, even when they act in weighted spaces.

We start by observing  that the  scale homogeneity of $S^{\alpha,\beta}$ is $2$ since a change of variables yields
\begin{align*}
S^{\alpha,\beta}(t)\left(I_s f\right)=I_s \left(S^{\alpha,\beta}(s^2t) f\right),\qquad I_sf(y)=f(sy),\qquad t, s>0,
\end{align*}
which in particular gives 
\begin{align*}
S^{\alpha,\beta}(t)f =I_{1/\sqrt t} \left(S^{\alpha,\beta}(1) I_{\sqrt t} f\right),\qquad t> 0.
\end{align*}
The boundedness of $S^{\alpha,\beta}(t) $ in $L^p_m$ is then equivalent to that for $t=1$ and $\|S^{\alpha,\beta}(t)\|_p= \|S^{\alpha,\beta}(1)\|_p$ and this is equivalent to 
$\alpha <\frac{M+m}p < M-\beta$ (in  particular $\alpha+\beta<M$), see Proposition \ref{Boundedness theta}, with $\theta=0$.

For future purpose we also  observe that the adjoint of $S^{\alpha,\beta}(t)$ taken with respect to the measure $\mu_m$
is given by  the operator 
\begin{align}\label{Adjoint integral operator}
\left(S^{\alpha,\beta}(t)\right)^{\ast m}f(y)=t^{-\frac {M} 2}\,\left (\frac{|y|}{\sqrt t}\wedge 1 \right)^{-\beta}\int_{\R^N}\left(\frac{|y|}{|z|}\right)^{-m}  \left (\frac{|z|}{\sqrt t}\wedge 1 \right)^{-\alpha}
\exp\left(-\frac{|y-z|^2}{\kappa t}\right)f(z) \,dz.
\end{align}

\subsection{$\mathcal{R}$-boundedness when $M+m>$0}

If   $d$ denotes the euclidean distance, $\left(\R^M,d,\mu_m\right)$ is  of homogeneous type.  In what follows we write $A_p(\mu_m)$, $RH_p(\mu_m)$, $\mathcal{M}_{\mu_m}$ to denote respectively the class of Muckenhoupt weights, the reverse H\"older class and the maximal function over balls taken with respect to the measure $\mu_m$. When $m=0$ we  write $A_p$, $RH_p$, $\mathcal{M}$.

We observe preliminarily that  \eqref{gen eq semigrop est Lp_m 1} yields 
when $m\geq 0$
\begin{align}\label{gen eq semigrop est Lp_m 2}
|S^{\alpha,\beta}(t)f(y)|&\leq C t^{-\frac {M+m} 2}\,\left (\frac{|y|}{\sqrt t}\wedge 1 \right)^{-\alpha}\int_{\R^M}  \left (\frac{|z|}{\sqrt t}\wedge 1 \right)^{-\beta-m}
\exp\left(-\frac{|y-z|^2}{\kappa t}\right)|f(z)| \,d\mu_m(z).
\end{align}
This follows after observing that, when $m\geq 0$, one has 
\begin{align*}
|z|^{-m}\left (|y|\wedge 1\right )^{-\beta}\leq \left (|z|\wedge 1\right )^{-\beta-m}, \quad z\in\R^M\setminus\{0\}.
\end{align*}
When  $m<0$, up to a small perturbation of the constant in the exponential argument,  estimate \eqref{gen eq semigrop est Lp_m 2} continues to hold in the range $\frac{|y|}{\sqrt t}\leq 1$, $\frac{|z|}{\sqrt t}\geq 1$. Indeed  in this case one has, for $\epsilon>0$ and for some $K>0$,
\begin{align*}
|z|^{-m}\exp\left(-\eps|y-z|^2\right)\leq K, \quad |y|\leq 1, \;|z|\geq 1.
\end{align*}
This implies, for $\frac{|y|}{\sqrt t}\leq 1$, $\frac{|z|}{\sqrt t}\geq 1$ and $\kappa'>\kappa$
\begin{align}\label{gen eq semigrop est Lp_m 3}
\nonumber |S^{\alpha,\beta}(t)f(y)|&\leq C t^{-\frac {N} 2}\,\left (\frac{|y|}{\sqrt t}\right)^{-\alpha}\int_{\R^M}  
\exp\left(-\frac{|y-z|^2}{\kappa t}\right)|f(z)| \,dz\\[1ex]
&\leq CK t^{-\frac {M+m} 2}\,\left (\frac{|y|}{\sqrt t}\right)^{-\alpha}\int_{\R^M}  
\exp\left(-\frac{|y-z|^2}{\kappa' t}\right)|f(z)| \,d\mu_m(z).
\end{align}

We  prove the $\mathcal{R}$-boundedness of the family $\left(S^{\alpha,\beta}(t)\right)_{t\geq 0} $ using the extrapolation result of Theorem \ref{Extrap Mack}. We follow the proof in \cite[Theorem 2.9]{bui} but new complications arise  because the operator is non-symmetric and the measure $\mu_m$ is not the Lebesgue one. In particular we have to distinguish between the cases $m \ge 0$ and $-M<m<0$ and both the maximal functions with respect to the Lebesgue measure and the weighted one appear. Note that,  since $M+m>0$,  condition (iii) in Proposition \ref{Boundedness theta} implies $\beta<M$ and $\alpha<M+m$.

 For the reader's convenience, in what follows we write for $t>0$, $B=B(0,\sqrt{t})$ and

\begin{align}\label{gen pezzi}
\nonumber S^{\alpha,\beta}(t)f&=\chi_{B^c}\left(S^{\alpha,\beta}(t)f\chi_{B^c}\right)+\chi_{B}\left(S^{\alpha,\beta}(t)(f\chi_{B})\right)+\chi_{B^c}\left(S^{\alpha,\beta}(t)f\chi_B\right)+\chi_{B}\left(S^{\alpha,\beta}(t)(f\chi_{B^c})\right)\\[1ex]&:=
S^{\alpha,\beta}_1(t)f+S^{\alpha,\beta}_2(t)f+S^{\alpha,\beta}_3(t)f+S^{\alpha,\beta}_4(t)f.
\end{align}

\medskip
 \begin{prop} \label{gen Teo Extra mu_m}
 Let $M+m>0$, $1<p< \infty$ and assume that $\alpha<\frac{M+m}p < M-\beta$. Let 
\begin{eqnarray} \label{notazione}
q_0&=&\frac{M+m}{\alpha}\ {\rm when}\  \alpha>0, \quad  q_0=\infty \ {\rm  when\ } \alpha \leq 0 \\
 \nonumber p_0&=& \left(\frac{M+m}{\beta+m}\right)'\ {\rm when}\  \beta+m>0, \quad  p_0=1 \ {\rm  when\ } \beta+m \leq 0 
\end{eqnarray}
so that $p_0<p<q_0$.
Then for every weight  
$$w\in A_{\frac{p}{p_0}}(\mu_m)\cap RH_{\left (\frac{q_0}{ p}\right )'}(\mu_m)$$
  there exist $C>0$ depending on  the   $A_{\frac{p}{p_0}}(\mu_m)$ and  $RH_{\left (\frac{q_0}{ p}\right )'}(\mu_m)$ constants of $w$ such that for every  $t\geq 0$ one has 
$$\|S^{\alpha,\beta}(t)f\|_{L^p(w)}\leq C\|f\|_{L^p(w)},\quad f\in L^p(\R^M,w d\mu_m)=:L^p(w).$$

Finally,  if in addition $p_0<2<q_0$ (i.e. $S^{\alpha,\beta}(1)$ is bounded on $L^2(\R^M,w d\mu_m)$) then the family $\left(S^{\alpha,\beta}(t)\right)_{t\geq 0}$ is $\mathcal{R}$-bounded on $L^p(\R^M,w d\mu_m)$.
\end{prop}
 
We split the proof in four  lemmas according to  \eqref{gen pezzi}.

 \begin{lem}\label{gen lemma T1}
The estimate of  Proposition \ref{gen Teo Extra mu_m} holds for $(S^{\alpha,\beta}_1(t))_{t\geq 0}$.
\end{lem}
{\sc{Proof.}} 
Assume first that $m\geq 0$. Then using \eqref{gen eq semigrop est Lp_m 2} and Lemma \ref{convolution radial} with $p(y)=|y|^m$  we get
\begin{align*}
|S^{\alpha,\beta}_1(t)f(y)|&\leq C t^{-\frac {M+m} 2}\int_{\R^M} 
\exp\left(-\frac{|y-z|^2}{\kappa t}\right)|f(z)| \,d\mu_m(z)\leq C\mathcal{M}_{\mu_m}f(y),
\end{align*}
The claim then  follows  since $\mathcal M_{\mu_m}$ is bounded on  $L^p(w)$.

When $-M<m<0$ we  use \eqref{gen eq semigrop est Lp_m 1} (and Lemma \ref{convolution radial} with respect to the Lebesgue measure) to  get
\begin{align*}
|S^{\alpha,\beta}_1(t)f(y)|&\leq C t^{-\frac {M} 2}\int_{\R^M} 
\exp\left(-\frac{|y-z|^2}{\kappa t}\right)|f(z)|\chi_{B^c} (z)\,dy\leq C\mathcal{M}f(y),
\end{align*}
Since $w\in A_{\frac{p}{p_0}}(\mu_m)$, by Proposition \ref{Prop Macke} there exists   $r$ sufficiently close to $p_0$ such that $p_0<r<p<q_0$ and  $w\in A_{\frac{p}{r}}(\mu_m)$.  Since $-M<m<0$, Lemma \ref{radial weights} (i) gives $|y|^m\in A_r(dy)$ and  then  Lemma \ref{Max functs estimate} yields
\begin{align*}
|S^{\alpha,\beta}_1(t)f(y)|&\leq C\left(\mathcal{M}_{\mu_m}|f|^{r}(y)\right)^{\frac 1 r }.
\end{align*}
Since $w\in A_{\frac p r}(\mu_m)$, $\mathcal M_{\mu_m}$ is bounded on $L^{\frac p r}(w)$ and we get $\|S^{\alpha,\beta}_1(t)f\|_{L^p(w)}\leq C\|f\|_{L^p(w)}$.\qed

  \begin{lem}\label{gen lemma T2}
The estimate of  Proposition \ref{gen Teo Extra mu_m} holds for $(S^{\alpha,\beta}_2(t))_{t\geq 0}$.
 \end{lem}
{\sc{Proof.}}  
Using \eqref{gen eq semigrop est Lp_m 1} 
and H\"older's inequality we get 
\begin{align*}
|S^{\alpha,\beta}_2(t)f(y)|&\leq C t^{-\frac {M+m} 2}\left (\frac{|y|}{\sqrt t}\right)^{-\alpha}\int_{B}  \left (\frac{|z|}{\sqrt t}\right)^{-\beta-m}
|f(z)| \,d\mu_m(z) \\
&\leq Ct^{-\frac{M+m}{2}}\left (\frac{|y|}{\sqrt t}\right)^{-\alpha}\|f\|_{L^p(w)} \left(\int_{B}\left (\frac{|z|}{\sqrt t} \right)^{-(\beta+m)p'}w(z)^{1-p'}\ d\mu_m(z)\right)^{\frac{1}{p'}}.
\end{align*}

Setting $v=w^{1-p'}$ this implies
\begin{align*}
\|S^{\alpha,\beta}_2(t)f\|_{L^p(w)}^p\leq Ct^{-\frac{M+m}{2}p}\|f\|_{L^p(w)}^p \left(\int_{B}\left (\frac{|z|}{\sqrt t} \right)^{-(\beta+m)p'}v(z)\ d\mu_m(z)\right)^{\frac{p}{p'}} \int_{B} \left (\frac{|y|}{\sqrt t}\right)^{-\alpha p}w(y)\ d\mu_m(y).
\end{align*}
Let us treat  the first integral.  If $\beta+m>0$, then  one has
\begin{align*}
\int_{B}\left (\frac{|z|}{\sqrt t} \right)^{-(\beta+m)p'}v(z)\ d\mu_m(z)&=\sum_{j\geq 0}\int_{2^{-j-1}\leq\frac{|z|}{\sqrt t}<2^{-j}}\left (\frac{|z|}{\sqrt t} \right)^{-(\beta+m)p'}v(z)\ d\mu_m(z)
\\[1ex]
&\leq C \sum_{j\geq 0} 2^{j(\beta+m)p'}v(2^{-j}B).
\end{align*}
By property (vii) of Proposition \ref{Prop Macke}, $v\in A_{\frac{p'}{q'_0}}\cap RH_{\left (\frac{p'_0}{p'}\right)'}$;  by property (v) of Proposition \ref{Prop Macke}  there exists $r>p'$ such that $v\in RH_{\left (\frac{p'_0}{r}\right)'}$.
Lemma \ref{Doubling weight} then implies
$$v(2^{-j}B)\leq C v(B)\left(\frac{\mu_m\left(2^{-j}B\right)}{\mu_m\left(B\right)}\right)^{\frac{(\beta+m) r}{M+m}}=C  v(B)2^{-jr(\beta+m)}.$$
Therefore since $\beta+m> 0$
\begin{align*}
\int_{B}\left (\frac{|z|}{\sqrt t} \right)^{-(\beta+m)p'}v(z)\ d\mu_m(z)&\leq Cv(B) \sum_{j\geq 0}   2^{-j(\beta+m)(r-p')}=Cv(B).
\end{align*}
The last inequality holds also  when $\beta+m\leq 0$, since in this case  
\begin{align*}
\int_{B}\left (\frac{|z|}{\sqrt t} \right)^{-(\beta+m)p'}v(z)\ d\mu_m(z)\leq \int_{B}v(z)\ d\mu_m(y)=v(B).
\end{align*}
Similarly if $\alpha>0$ then 
\begin{align*}
\int_{B} \left (\frac{|y|}{\sqrt t}\right)^{-\alpha p}w(y)\ d\mu_m(y)&=\sum_{j\geq 0}\int_{2^{-j-1}\leq\frac{|y|}{\sqrt t}<2^{-j}}\left (\frac{|y|}{\sqrt t}\right)^{-\alpha p}w(y)\ d\mu_m(y)\\[1ex]
&\leq C \sum_{j\geq 0} 2^{j\alpha p}w(2^{-j}B).
\end{align*}
Since $w \in A_{\frac{p}{p_0}}\cap RH_{\left (\frac{q_0}{ p}\right )'}$  by property (v) of Proposition \ref{Prop Macke}  there exists $r>p$ such that $w\in  RH_{\left (\frac{q_0}{ r}\right )'}$.
By Lemma \ref{Doubling weight} then 
$$w(2^{-j}B)\leq C w(B)2^{-jr\alpha}.$$
Therefore
\begin{align*}
\int_{B}\left (\frac{|y|}{\sqrt t} \right)^{-\alpha p}w(y)\ d\mu_m(y)\leq C w(B) \sum_{j\geq 0}  2^{-j\alpha(r-p)}=Cw(B).
\end{align*}
The last inequality holds also  when $\alpha\leq 0$, since in this case  
\begin{align*}
\int_{B}\left (\frac{|y|}{\sqrt t} \right)^{-\alpha}w(y)\ d\mu_m(y)\leq \int_{B}w(y)\ d\mu_m(y)=w(B).
\end{align*}
Putting together the last inequalities we have in any case
$$\|S^{\alpha,\beta}_2(t)f\|_{L^p(w)}^p\leq C \|f\|_{L^p(w)}^p t^{-p\frac{M+m}{2}}\left (v(B)\right)^{\frac{p}{p'}}w(B).$$
Since $\beta<M$ from property (i) of Proposition \ref{Prop Macke} we get  $w\in A_{\frac{p}{p_0}}\subseteq A_p$  which implies, by the definition \eqref{Def A_p} of $A_p$ weights,
$\sup_{ t >0}t^{-p\frac{M+m}{2}}\left (v(B)\right)^{\frac{p}{p'}}w(B) <\infty$.
\qed
 
\begin{lem}\label{gen lemma T3}
The estimate of  Proposition \ref{gen Teo Extra mu_m} holds for $(S^{\alpha,\beta}_3(t))_{t\geq 0}$.
 \end{lem}
{\sc{Proof.}}
Using \eqref{gen eq semigrop est Lp_m 1} we get
\begin{align*}
|S^{\alpha,\beta}_3(t)f(y)|&\leq C t^{-\frac {M+m} 2}\int_{B}  \left (\frac{|z|}{\sqrt t} \right)^{-\beta-m}
\exp\left(-\frac{|y-z|^2}{\kappa t}\right)|f(z)| \,d\mu_m(z).
\end{align*}
Let us fix $r$ such that $p_0'<r<p<q_0$.
Applying H\"older's inequality we obtain
\begin{align*}
|S^{\alpha,\beta}_3(t)f(y)|&\leq C  \left(t^{-\frac{M+m}{2}} \int_{\R^M}
\exp\left(-\frac{|y-z|^2}{\kappa t}\right)|f(z)|^r\ d\mu_m(z)\right)^\frac{1}{r} \\[1ex]
&\hspace{20ex}\times\left(t^{-\frac{M+m}{2}} \int_{B}\left (\frac{|z|}{t^{\frac{1}{2}}} \right)^{-(\beta+m)r'}
\ d\mu_m(z)\right)^\frac{1}{r'}.
\end{align*}
The substitution $\xi=z/\sqrt t$ and Lemma \ref{convolution radial} yield 
\begin{align*}
|S^{\alpha,\beta}_3(t)f(y)|&\leq C  \left(\mathcal{M}_{\mu_m}|f|^r(y)\right)^{\frac 1 r} \left(\int_{B(0,1)}|\xi|^{-(\beta+m)r'+m}\ d\xi\right)^\frac{1}{r'}=C  \left(\mathcal{M}_{\mu_m}|f|^r(y)\right)^{\frac 1 r}
\end{align*}
Since $w\in A_{\frac{p}{p_0}}$, by Proposition \ref{Prop Macke} there exists   $r$ sufficiently close to $p_0$ such that $p_0<r<p<q_0$ and  $w\in A_{\frac{p}{r}}$ . This implies that  $\mathcal M_{\mu_m}$ is bounded on $L^{\frac p r}(w)$ which, using the latter inequality, proves the required claim.\\\qed
 \medskip
 
 Finally, in order to prove the boundedness of $S^{\alpha,\beta}_4(t)$, we employ estimates \eqref{gen eq semigrop est Lp_m 2}, \eqref{gen eq semigrop est Lp_m 3} which allow to  equivalently prove, up to a small modification of the constant in the exponential argument,  the  boundedness of the operator
\begin{align*}
 F_4(t)f(y)
&=\chi_{B}(y)\,{t}^{-\frac {M+m} 2}\left(\frac{|y|}{\sqrt{t}}\right)^{-\alpha}\,\int_{B^c} \exp\left(-\frac{|y-z|^2}{\kappa t}\right)|f(z)| \,d\mu_m(z).
\end{align*}

We  apply a duality argument and we observe that the adjoint of $F_4(t)$ taken with respect to the measure $\mu_m$
is given by  the operator
\begin{align*}
\left(F_4(t)^{\ast m}\right)f(y)&= t^{-\frac {M+m} 2}\int_{B}  \left (\frac{|z|}{\sqrt t} \right)^{-\alpha}
\exp\left(-\frac{|y-z|^2}{\kappa t}\right)|f(z)| \,d\mu_m(z)\\[1ex]
&=S^{\beta,\alpha}_3(t)f(y)
\end{align*}
 With this aim let us observe that $q_0'<p'<p_0'$.

 \begin{lem}\label{gen lemma T4}
The estimate of  Proposition \ref{gen Teo Extra mu_m} holds for   $(S^{\alpha,\beta}_4(t))_{t\geq 0}$.
 \end{lem}
{\sc{Proof.}} 
We apply a duality argument. Let $g\in L^{p'}(\R^M,w\mu_m)$; since $(F_4(t))^{\ast m}=S^{\beta,\alpha}_3(t)$ we obtain
\begin{align*}
\int_{\R^M}F_4(t)fg\, w\mu_m=\int_{\R^M}f\,S^{\beta,\alpha}_3(t)(g w)\mu_m=\int_{\R^M}f\,\frac{S^{\beta,\alpha}_3(t)(g w)}{w}w\mu_m.
\end{align*}
Using H\"older's inequality we then yield
\begin{align*}
\left|\int_{\R^M}F_4(t)fg\, w\mu_m\right|\leq\|f\|_{L^p(\omega)}\left\|\frac{S^{\beta,\alpha}_3(t)(g w)}{w}\right\|_{L^{p'}(\omega)}=\|f\|_{L^p(\omega)}\left\|S^{\beta,\alpha}_3(t)(g w)\right\|_{L^{p'}(\omega^{1-p'})}.
\end{align*}
By property (vii) of Proposition \ref{Prop Macke}, $\omega^{1-p'}\in A_{\frac{p'}{q_0'}}\cap RH_{\left(\frac{p_0'}{ p'}\right)'}$. Then using the estimate for $S^{\beta,\alpha}_3$ and with $p$ replaced by $p'$ we get
\begin{align*}
\left|\int_{\R^M}F_4(t)fg\, w\mu_m\right|\leq C\|f\|_{L^p(\omega)}\left\|g w\right\|_{L^{p'}(\omega^{1-p'})}=C\|f\|_{L^p(\omega)}\left\|g \right\|_{L^{p'}(\omega)}
\end{align*}
which concludes the proof.
\qed

 \medskip
 We can finally prove Proposition \ref{gen Teo Extra mu_m}.
 \medskip
 
 {\sc \bf (Proof of Proposition \ref{gen Teo Extra mu_m})} The first claim follows by using \eqref{gen pezzi} and Lemmas \ref{gen lemma T1}, \ref{gen lemma T2}, \ref{gen lemma T3}, \ref{gen lemma T4}. If  $p=2$ satisfies $p_0<2<q_0$, then the  $\mathcal{R}$-boundedness of $(S^{\alpha,\beta}(t))_{t\geq 0}$  on $L^p(\R^M,w d\mu_m)$ follows by Theorems \ref{Square funct R-bound}, \ref{Extrap Mack} since, in this case, the boundedness  of $S^{\alpha,\beta}(t)$ in all the Muckenhoupt weighted spaces just proved implies the square function estimate  required by Theorem \ref{Square funct R-bound}.\qed

\subsection{$\mathcal{R}$-boundedness in the general case}

Let us remove the assumptions $M+m>0$ and  $p_0<2<q_0$ first in  the case of the Lebesgue measure, that is when $m=0$. Since we use different measures here we do not shorten $L^p(\R^M, |y|^mdy)$ to $L^p_m$.

 \begin{teo} \label{gen R-bound sempre}
 Let $1<p< \infty$ and let us suppose that $\alpha<\frac{M+m}p < M-\beta$.  Then$\left(S^{\alpha,\beta}(t)\right)_{t\geq 0}$ is $\mathcal{R}$-bounded on $L^p(\R^M)$.

\end{teo}
{\sc{Proof.}} 
If $\alpha<\frac M 2<M-\beta$ i.e. $\alpha,\beta<\frac M 2$, the thesis is part of   Proposition  \ref{gen Teo Extra mu_m} with $m=0$. Let us suppose now $\alpha\geq \frac M2$ or $\beta\geq \frac M2$; in this case
$S^{\alpha,\beta}(t)$ is not bounded on $L^2\left(\R^M\right)$ and therefore $p\neq 2$. 

Given $m\in\R$,   let us consider the isometry 
\begin{align*}
T_{\frac m p}: L^p(\R^M,|y|^m dy)\to L^p(\R^M,dy),\quad f\mapsto |y|^{\frac m p}f.
 \end{align*}
A straightforward computation shows that
\begin{align*}
T_{-\frac m p} S^{\alpha,\beta}(t)T_{\frac m p}f=\tilde S^{\alpha,\beta,\frac{m}p}(t)f
\end{align*}
where $\tilde S^{\alpha,\beta,\frac{m}p}$ is the operator defined in  \eqref{Equiv integral oper} with $r=m/p$. Lemma \ref{Lem Equiv integral oper} gives
\begin{align}\label{equiv gen m}
|T_{-\frac m p} S^{\alpha,\beta}(t)T_{\frac m p}f|\leq C S^{\alpha+\frac m p,\beta-\frac{m}p}(t)|f|
\end{align}

 Since  the $\mathcal{R}$-boundedness of a family of operators is preserved by isometries and pointwise domination, from the equality \eqref{equiv gen m} one can easily deduce that  $\left(S^{\alpha,\beta}(t)\right)_{t\geq 0}$ is $\mathcal{R}$-bounded on $L^p(\R^M)$ if there exists $m\in\R$ such that $\left(S^{\alpha+\frac m p,\beta-\frac{m}p}(t)\right)_{t\geq 0}$ is $\mathcal{R}$-bounded on $L^p(\R^M, |y|^mdy)$.
From Proposition \ref{gen Teo Extra mu_m} it  is then sufficient to require 
\begin{align*}
 0&<M+m,\qquad\alpha+\frac mp <\frac{M+m}{2}< M-\beta+\frac m p.
 \end{align*} 
%
%
By elementary calculation the latter inequalities read as
\begin{align*}
\begin{cases}
M+m>0;\\[1.3ex]
m\left(\frac 1 p-\frac 1 2\right)<\frac M 2-\alpha;\\[1.3ex]
m\left(\frac 1 p-\frac 1 2\right)>\beta-\frac M 2,
\end{cases}
\end{align*}
 If $p<2$ the  system has a solution $m$ when 
$$\beta-\frac M 2<-M\left ({\frac{1}{p} -\frac12}\right ) <\frac{M}{2} -\alpha
$$
that is when $\alpha <\frac{M}{p} <M-\beta$.  If $p>2$  the  claim follows in the same way.
%
\qed

 The results for $S^{\alpha,\beta}(t)$ in $L^p\left(\R^M,d\mu_m\right)$ are immediate consequence of those of $S^{\alpha-\frac m p,\beta+\frac mp}(t)$ in $L^p(\R^M, dy)$.  Note that the condition $M+m>0$ is no longer required.

 
\begin{teo} \label{gen R-bound sempre m}
 If , $1<p< \infty$ and $\alpha<\frac{M+m}p < M-\beta$, then the family
 $\left(S^{\alpha,\beta}(t)\right)_{t\geq 0}$ is $\mathcal{R}$-bounded on $L^p_m=L^p(\R^M, d\mu_m)$.

\end{teo} {\sc{Proof.}} 
Let us consider the isometry 
\begin{align*}
T_{-\frac m p}: L^p(\R^M, dy)\to L^p(\R^M,|y|^mdy),\quad f\mapsto |y|^{-\frac m p}f.
 \end{align*}
Then, as done in the previous proof, one has,  using Lemma \ref{Lem Equiv integral oper},
\begin{align*}
|T_{\frac m p} S^{\alpha,\beta}(t)T_{-\frac m p}f|=|T^{\alpha,\beta,-\frac{m}p}f|\leq  C S^{\alpha-\frac m p,\beta+\frac{m}p}|f|.
\end{align*}
By construction,   the boundedness conditions for $S^{\alpha-\frac m p,\beta+\frac{m}p}(1)$ in $L^p(\R^M, dy)$ are  satisfied under the given hypotheses  on $S^{\alpha,\beta}$. Then the family  $S^{\alpha-\frac m p,\beta+\frac{m}p}$ 
 is $\mathcal{R}$-bounded  in $L^p(\R^M, dy)$  by Theorem \ref{gen R-bound sempre}; the same result then also holds for  $T_{\frac m p} S^{\alpha,\beta}T_{-\frac m p}$ by domination. By similarity this yields the  $\mathcal{R}$-boundedness of $S^{\alpha,\beta}$  in $L^p(\R^M, |y|^mdy)$.\qed

\section{Domain and maximal regularity for $\mathcal L=\Delta_x+L_y$}

\subsection{Basic facts}
Here we deduce generation results for the whole operator $\mathcal L=\Delta_x+L_y$ by standard tensor product arguments.
If $X,Y$ are function spaces over $G_1, G_2$ we denote by $X\otimes Y$ the algebraic tensor product of $X,Y$, that is the set of all functions $u(x,y)=\sum_{i=1}^n f_i(x)g_i(y)$ where $f_i \in X, g_i \in Y$ and $x \in G_1, y\in G_2$.
If $T,S$ are linear operators  on $X,Y$ we denote by $T\otimes S$ the operator on $X\otimes Y$ defined by 
$$
 \left (\left (T\otimes S  \right  )u\right)(x,y)=\sum_{i=1}^n (T f_i)(x)(Sg_i)(y)
$$
and we keep the same notation to denote its bounded extension to the completion of $X\otimes Y$, if no ambiguity can arise.
The generation result for $\mathcal L$ follows from well-known general facts.
We start by two preliminary lemmas where the tensor product of a semigroup $(e^{tA})_{t \ge0}$ with the identity operator is considered. The first follows from \cite[AI, Section 3.7]{nagel}.
\begin{lem} 
Let $(e^{tA})_{t \ge0}$ be the semigroup  generated by $(A,D(A))$ in $L^p(\Omega, d\mu)$. The family $(e^{tA})_{t \ge0}\otimes I)_{t \ge 0}$ on $L^p(\Omega\times\Lambda, d\mu \otimes d\nu)=\ov{L^p(\Omega, d\mu)\otimes L^p(\Lambda, d\nu)}$ is a semigroup  generated by the closure $\ov{A\otimes I}$ of the operator $A\otimes I$ initially defined on  $D(A)\otimes L^p(\Lambda)$.
\end{lem}
Let us introduce the operator $(A^\otimes, D(A^\otimes))$ 
\begin{align*}
D(A^\otimes)&:=\{u\in L^p(\Omega\times\Lambda):\ u(\,\cdot\,,y)\in D(A)\ \textrm{for almost every}\ y\in \Lambda,\ Au(\cdot,y)\in L^p(\Omega\times\Lambda)\}\\[1.5ex]
A^{\otimes}u(\,\cdot\,,y)&:=Au(\,\cdot\,,y),\quad \textrm{for almost every}\ y\in \Lambda.
\end{align*}

We can identify the domain of the generator as follows.
\begin{lem}
Keeping the notation of the previous Lemma, we have $\ov{A\otimes I}=A^\otimes$.
\end{lem}
{\sc Proof.} We start by proving that $A^\otimes$ is  a closed operator. Let $(u_n)\in D(A^\otimes)$ such that $u_n\to u$, $A^\otimes u_n\to v$ in $L^p(\Omega\times \Lambda)$. Up to considering a subsequence,
$u_n(\cdot,y)\to u(\cdot,y)$, $Au_n(\cdot,y)\to v(\cdot,y)$ for almost  every $y\in\Lambda$. Then, since $A$ is closed, $u(\cdot,y)\in D(A)$ and $v(\cdot,y)=Au(\cdot,y)$ for almost every $y\in \Lambda$. It follows $u\in D(A^\otimes)$ and $A^\otimes u=v$.
Next we prove that $\ov{A\otimes I}\subset A^\otimes$. By the definition, it easily follows that $A\otimes I\subset A^\otimes$. Since $A^\otimes$ is closed by the previous step, the inclusion  $\ov{A\otimes I}\subset A^\otimes$ is proved. Finally we prove that, for $\lambda$ large enough, the operator $\lambda-A^\otimes$ is injective. Indeed, if $u\in D(A^\otimes)$ and $\lambda u-A^\otimes u=0$, then $\lambda u(\cdot, y)- Au(\cdot, y)=0$ for almost every $y\in\Lambda$ and, by the injectivity of $A$, $u(\cdot,y)=0$ for almost every $y\in\Lambda$.\qed

We therefore deduce the following generation result, under the assumption that $A_{m,p}$ generates in $L^p_m(0,\infty)$. Here $A_{m,p}$ denotes the degenerate operator $L_{m,p}$ of Section 4,  with  Dirichlet boundary conditions, or the Bessel operator $B^n_{m,p}$ of Section 3, under Neumann boundary conditions. We still write $A_y$ to indicate that $A$ acts only in the $y$ variable. 
\begin{lem}
Let $(e^{z \Delta_x})_{z \in \Sigma_\phi}$, $(e^{z A_y})_{z \in \Sigma_\phi}$ be the semigroups  generated  by $\left(\Delta_x,W^{2,p}(\R^N)\right)$ in $L^p(\R^N)$, $\left(A_{m,p},D(A_{m,p})\right)$ in $L_m^p(0,\infty)$, respectively.  The operators $\Delta_x^\otimes$  and $A_{m,p}^\otimes$ defined by
\begin{align*}
D(\Delta_x^\otimes):=\Big\{&u\in L^p_m(\R_+^{N+1}):\ u(\,\cdot\,,y)\in W^{2,p}(\R^N)\ \textrm{for a.e.}\ y\in (0,\infty),\ \nabla_xu(\cdot,y),\Big.  \\ 
& \Big.D^2_x u (\cdot,y)\in L^p_m(\R_+^{N+1}),\quad  \Delta_x^{\otimes}u(\,\cdot\,,y):=\Delta_x u(\,\cdot\,,y),\quad \textrm{for almost every}\ y\in (0,\infty)\Big\};
\end{align*}
\begin{align*}
D(A_{m,p}^\otimes):=\Big\{&u\in L^p_m(\R^{N+1}_+):\ u(x,\cdot)\in D(A_{m,p})\ \textrm{for a.e.}\ x\in \R^N,\ A_yu(x,\cdot)\Big.\\\Big.&\in L^p_m(\R_+^{N+1}), \quad A_{m,p}^{\otimes}u(x, \cdot):=A_yu(x, \cdot),\quad \textrm{for almost every}\ x\in \R^N\Big\}
\end{align*}
generate the semigroups $(e^{z \Delta_x} \otimes I)_{z \in \Sigma_\phi}$, $( I\otimes e^{z A_y})_{z \in \Sigma_\phi}$ in  $L^p_m(\R_+^{N+1})$.
\end{lem}

\subsection{Maximal regularity and domain characterization }
We can finally prove maximal regularity and domain characterization for $\mathcal L=\Delta_x+L_y$ and $\mathcal L=\Delta_x+B^n_y$. Both cases have similar proofs but some details are different since the domain of $B^n$ is more regular. In the gradient estimates for $B^n$, in fact, the factor  $y^{-s_1-1}$ does not appear, see Proposition \ref{Estimates gradient kernel bessel}, in contrast with Proposition \ref{kernelL} where it is assumed that $s_1 \neq 0$. However, if $s_1=0$, then $b=0$ and $c \geq 1$ so that $L=B^d=B^n$. Therefore, we distinguish between the cases of $L$ with $s_1 \neq 0$ and $B^n$, the case of $L$ with $s_1=0$ being included in the last.

\subsection{$ \mathcal L=\Delta_x+L_y$ with $s_1 \neq 0$}
First we state a $\mathcal R$-boundedness result for the heat kernel of $L$ and its gradient. We remark that the $\mathcal R$-boundedness of the heat kernel has been also proved in \cite{met-negro-spina 7}.
\begin{teo} \label{Rbounded} If $s_1<\frac{m+1}{p} <s_2+2$, then the  family $(e^{zL_{m,p}})_{z \in \Sigma_\phi}$, $\phi<\pi/2$ is $\mathcal R$-bounded in $L^p_m(0,\infty)$,  hence $L_{m,p}$ has maximal regularity.

If $s_1+1<\frac{m+1}{p}<s_2+2$, then families  $(\frac{\sqrt z}{y} e^{zL_{m,p}})_{z \in \Sigma_\phi}$, $(\sqrt z D_ye^{zL_{m,p}})_{z \in \Sigma_\phi}$, $\phi<\pi/2$ are $\mathcal R$-bounded in $L^p_m(0,\infty)$.
\end{teo}
{\sc Proof.} By Proposition \ref{kernelL} we have
$$
|e^{zL_{m,p}}f| \leq C S^{\alpha, \beta}(c|z|)|f|
$$ pointwise, for $\alpha=s_1$, $\beta=s_1-c$  and suitable positive constants $C ,c$ (note also that since $s_1+s_2=c-1$, then $1-\beta=s_2+2$).  
The assertion for   $(e^{zL_{m,p}})_{z \in \Sigma_\phi}$ then  follow from Theorem \ref{gen R-bound sempre m} together with Corollary \ref{domination}. 

Those for $(\frac{\sqrt z}{y} e^{zL_{m,p}})_{z \in \Sigma_\phi}$, $(\sqrt z D_ye^{zL_{m,p}})_{z \in \Sigma_\phi}$ are proved in a similar way, using Proposition \ref{kernelL} and setting
setting $\alpha=s_1+1$ and $\beta=s_1-c$.  
\qed

\begin{prop} \label{analyt} Let $s_1<\frac{m+1}{p}< s_2+2$.
Then the closure of the operator $\mathcal L$, initially defined on $W^{2,p}(\R^N)\otimes  D (L_{m,p})$,  generates a bounded analytic semigroup of angle $\pi/2$, $(e^{z \mathcal L})_{ z \in C_+}$, in $L^p_m(\R_+^{N+1})$ which has maximal regularity.

\end{prop}
{\sc Proof.} Observe first that 
$$\mathcal L=\Delta_x\otimes I +I\otimes L_y$$ on $W^{2,p}(\R^N)\otimes  D(L_{m,p})$. The family $(e^{z\Delta_x} \otimes e^{zL_y})_{z \in \C_+}$ is a semigroup and leaves $W^{2,p}(\R^N)\otimes  D(L_{m,p})$ invariant. This last, being dense, is then a core for the generator. The $\mathcal R$-boundedness of the family $(e^{z \mathcal L})_{z \in \Sigma_\phi}=(e^{z\Delta_x} \otimes e^{zL_y})_{z \in \Sigma_\phi}$, $\phi<\pi/2$ follows by composition  writing $(e^{z\Delta_x} \otimes e^{zL_y})= (e^{z\Delta_x} \otimes I)\circ (I \otimes e^{zL_y})$, 
using the above theorem and the $\mathcal R$-boundedness of $(e^{z \Delta_x})_{z \in \Sigma_\phi}$.
\qed
We 
note that, by construction, $e^{t\mathcal L}$ consists of integral operators. For $t>0$, $z_1=(x_1,y_1),z_2=(x_2,y_2)\in \R^{N+1}_+$

\begin{align*}
\nonumber e^{t\mathcal L}f(z_1)&=\int_{R^{N+1}_+}p(t,z_1,z_2)f(z_2)d m(z_2),\quad f\in L^p_m\left(\R^{N+1} \right ) \\[2ex]
p(t,z_1,z_2)&=(4\pi t)^{-\frac N 2}e^{-\frac{|x_1-x_2|^2}{4t}}p_{L_y}(t,y_1,y_2) \\
&\simeq t^{-\frac{N+1}{2}}\left (\frac{|y_1|}{t^{\frac{1}{2}}}\wedge 1 \right)^{-s_1} \left (\frac{|y_2|}{t^{\frac{1}{2}}}\wedge 1 \right)^{-s_1+c}
\exp\left(-\frac{|z_1-z_2|^2}{\kappa t}\right).
\end{align*}

Before describing the domain  of $\mathcal L_{m,p}$, let us show how the results for the $1d$ operator $L_y$ give easily a core.
\begin{prop} \label{core1}
 If $s_1<\frac{m+1}{p} < s_2+2$, then $$\mathcal D=\left \{u=y^{-s_1}v: v \in C_c^\infty (\R^{N} \times [0, \infty)), \ D_y v(x,0)=0 \right \}$$ is a core for $D(\mathcal L_{m,p})$.
\end{prop}
{\sc Proof.} $C_c^\infty (\R^N)$ is a core for $\Delta_x$ and $\mathcal D_1= \left \{u=y^{-s_1}v: v \in C_c^\infty [0,\infty), \ D_y v(0)=0 \right \}$ is a core for $L_{m,p}$, by Proposition \ref{core}. Then $C_c^\infty (\R^N)\otimes \mathcal D_1 \subset \mathcal D$ is dense in $W^{2,p}(\R^N)\otimes D(L_{m,p})$ for the  norm $\|u\|=\|u\|_{L^p_m}+\|\Delta_x u\|_{L^p_m}+\|L_y u\|_{L^p_m}$, hence for the graph norm induced by $\mathcal L$. Since  $W^{2,p}(\R^N)\otimes D(L_{m,p})$ is a core for $\mathcal L_{m,p}$, the proof is complete.
\qed


\begin{teo}  \label{Domain}
Let  $D\geq 0$, $s_1<\frac{m+1}{p}<s_2+2$. Then 
$$D(\mathcal L_{m,p})=\Big\{u \in W^{2,p}_{loc}(\R^{N+1}_+): u, \nabla_x u, D^2_x u, L_y u \in L^p_m(\R^{N+1}_+)\Big\}.
$$
\end{teo}
{\sc Proof. } Observe that, by construction, 
\begin{equation} \label{inclusione}
\mathcal{S}(\R^N)\otimes  D (L_{m,p})\subset W^{2,p}(\R^N)\otimes  D (L_{m,p}) \subset D(\Delta_x^\otimes) \cap D(L_{m,p}^\otimes) \subset D(\mathcal L_{m,p})
\end{equation}
where $\mathcal{S}(\R^N)$ denotes the Schwartz class.  Note that  $\mathcal S(\R^N)\otimes  D(L_{m,p})$ is a core for $\mathcal L_{m,p}$ by the above proposition (or also since it is invariant for $(e^{z\Delta_x} \otimes e^{zL_y})_{z \in \C_+}$).

 We endow  $D(\mathcal L_{m,p})$ with the graph norm and $Z:=D(\Delta_x^\otimes) \cap D(L_{m,p}^\otimes)$ with the norm $$\|u\|_Z=\|u\|_{L^p_m}+\|\Delta_x u\|_{L^p_m}+\|L_y u\|_{L^p_m},\ \ u \in Z,$$ so that the embedding $Z \subset D(\mathcal L_{m,p})$ is continuous.
Let us show that the graph norm and the norm of $Z$ are equivalent on $\mathcal{S}(\R^N)\otimes  D (L_{m,p})$.  Let $u \in \mathcal S(\R^N)\otimes  D (L_{m,p})$ and $f=\mathcal Lu$.     By taking the Fourier transform with respect to $x$ (with co-variable $\xi$) we obtain
$$
(-|\xi|^2+L_y)\hat u(\xi,\cdot)=\hat f(\xi,\cdot), \qquad |\xi|^2 \hat u(\xi, \cdot)=-|\xi|^2(|\xi|^2-L_y)^{-1}\hat f(\xi,\cdot).
$$
This means $\Delta_x u=-{\cal F}^{-1} M(\xi) {\cal F} f$, where ${\cal F}$ denotes the Fourier transform and $M(\xi)=|\xi|^2(|\xi|^2-L_y)^{-1}$.

The estimate $\|\Delta_x u\|_p \le C\|f\|_p$ (norms on $L^p_m(\R^{N+1})$) follows  from the boundedness of the multiplier $M$ in $L^p(\R^N; L^p_m(0, \infty))$ which we prove using Theorem \ref{mikhlin}.
%
%
%
 In fact, since $(e^{tL_y})_{ t\ge 0}$ is $\mathcal R$-bounded by Theorem \ref{Rbounded}, then the family
$$
\Gamma (\lambda)=\lambda (\lambda-L_y)^{-1}=\int_0^\infty \lambda e^{-\lambda t}e^{tL_y}\, dt, \quad \lambda>0
$$
is $\mathcal R$-bounded by \cite[Corollary 2.14]{KW} and  satisfies Mikhlin condition in Theorem \ref{mikhlin} for every $N$, by the resolvent equation  (or arguing as in  Lemma \ref{h} below). The same then holds for $M(\xi)=\Gamma(|\xi|^2)$, as readily verified.

The estimate $\|L_y u\|_p \le C\|f\|_p$ follows by difference and shows the equivalence of the graph norm and of the norm of $Z$ on $\mathcal{S}(\R^N)\otimes  D (L_{m,p})$. If $u \in D(\mathcal {L}_{m,p})$, let $(u_n) \subset \mathcal{S}(\R^N)\otimes  D (L_{m,p})$ converge to $u$ with respect to the graph norm. Then $(u_n)$ is a Cauchy sequence with respect to the norm of $Z$, hence  $u \in Z$ and the equivalence of the corresponding norms extends to $Z=D(\mathcal{L}_{m,p})$.
\qed
 A more detailed description of the domain of $\mathcal L$ follows immediately from the above theorem and the domain description of $D(L_{m,p})$ of Section 4. We do not list all the results that can be obtained in this way, since this is straightforward. See however Corollary \ref{Dreg} below for an important case. 

When $(m+1)/p >s_1+1$ and $u \in D(\mathcal L_{m,p})$,  the mixed derivatives $D_y\nabla_x u$ belong to $L^p_m$ even though $D_{yy} u$ could be not  in $L^p_m$. 

\begin{teo}   \label{Lp-estimates} Let  $D\geq 0$, $s_1+1<\frac{m+1}{p} < s_2+2$. Then 
$$D(\mathcal L_{m,p})=\{u \in W^{2,p}_{loc}(\R^{N+1}_+): u, \nabla_x u, D^2_x u, \frac{\nabla_x u}{y}, D_y\nabla_x u,  L_y u \in L^p_m(\R^{N+1}_+)\}.
$$
\end{teo}
{\sc Proof.} We proceed as in Theorem \ref{Domain} to estimate $D_y\nabla_x u$ for  $u \in \mathcal{S}(\R^N)\otimes  D (L_{m,p})$. We have
$$
(-|\xi|^2+L_y)\hat u(\xi,\cdot)=\hat f(\xi,\cdot), \qquad \xi D_y \hat u(\xi, \cdot)=-\xi D_y(|\xi|^2-L_y)^{-1}\hat f(\xi,\cdot).
$$ and this time the multiplier is
$$M(\xi)=-\xi D_y(|\xi|^2-L_y)^{-1}=-\xi\int_0^\infty e^{-|\xi|^2 t}D_y e^{tL_y}\, dt=-\xi\int_0^\infty \frac{e^{-|\xi|^2 t}}{\sqrt t}S_t\, dt
$$
where $(S_t)_{t \ge 0}=(\sqrt t D_y e^{tL_y})_{t \ge 0}$ is $R$-bounded, by Theorem \ref{Rbounded}. The Mikhlin condition of Theorem \ref{mikhlin} follows from \cite[Corollary 2.14]{KW} and the lemma below. 

The proof for $y^{-1} \nabla_x u$ is similar.
\qed

\begin{lem} \label{h}
Let $h(\xi,t)=\frac{\xi}{\sqrt t}e^{-|\xi|^2 t}$. Then if $|\alpha|=k$
$$|\xi|^k\int_0^\infty |D^\alpha_\xi h(\xi,t)|dt  \le C_k.
$$
\end{lem}
{\sc Proof. } Let $g(\eta)=\eta e^{-|\eta|^2}$ so that $h(\xi,t)=\frac{1}{t} g(\xi \sqrt t)$ 
. Let us observe that for $|\alpha|=k>0$ is $|D^\alpha g(\eta)|=P^{k+1}(\eta)e^{-|\eta|^2}$ where $P^{k+1}$ is a vector  polynomial of degree $k+1$. This implies in particular that
$$|\eta|^{k-1}|D^\alpha_\eta g(\eta)|\leq C_k  e^{-|\eta|^2},\qquad \text{for}\quad k\geq 0$$
which yields
$$|D^\alpha_\xi h(\xi,t)|=t^{\frac{k}{2}-1}| D^\alpha_\xi g(\xi \sqrt t)| \le C_k t^{-\frac{1}{2}}|\xi|^{1-k}e^{-|\xi|^2t}.$$
Then one has 
\begin{align*}
|\xi|^k\int_0^\infty |D^\alpha_\xi h(\xi,t)|dt  &\le C_k |\xi|\int_0^\infty t^{-\frac{1}{2}}e^{-|\xi|^2t}dt=C_k\int_0^\infty s^{-\frac 1 2 }e^{-s}ds= C_k\sqrt \pi.
\end{align*}
%
\qed

Finally, if $(m+1)/p>s_1+2$, then $D(\mathcal L_{m,p})$  has the maximal regularity one can expect.
\begin{cor} \label{Dreg}
Let  $D> 0$, $s_1+2<\frac{m+1}{p} <s_2+2$. Then 
$$D(\mathcal L_{m,p})=\{u \in W_m^{2,p}(\R^{N+1}_+): \frac{u}{y^2}, \frac{\nabla u}{y} \in L^p_m(\R^{N+1}_+)\}$$.
\end{cor}
{\sc Proof. } We have only to show that $y^{-2} u, y^{-1} D_y u, D_{yy} \in L^p_m$, since the rest follows from Theorem \ref{Lp-estimates}. Using Proposition \ref{gen-dom-subcritBDirichlet1} with $\theta=1$ we get
$$
\int_0^\infty \left (|D_{yy}u(x,y)|^p+\frac{|D_y u(x,y)|^p}{y^p}+\frac{|u(x,y)|^p}{y^{2p}} \right ) y^m dy \le C\int_0^\infty |L_y u(x,y)|^p y^m dy
$$
(the additional term containing $u$ on the right hand side does not appear, by homogeneity reasons). Integrating with respect to $x \in \R^N$ and using the estimate $\|L_y u\|_{L^p_m (\R^{N+1}_+)} \le C \|\mathcal{L} u\|_{L^p_m (\R^{N+1}_+)}$ the proof is complete.
\qed

\subsection{ $\mathcal L=\Delta_x+L_y$  with $s_1=0$ and $\mathcal L^n=\Delta_x+B^n_y$ }
If $s_1=0$ then $b=0$, $c \geq 1$ and   $L=B=D_{yy} +\frac{c}{y} D_y$ is a Bessel operator. Since $c \geq 1$, then $B^d=B^n$, see Section 2, and therefore it is sufficient to deal with $B^n$. Note, however, that $s_1=c-1$ for $B^n$ when $c<1$.

\begin{teo} \label{Rbounded1} If $c>-1$ and  $\frac{m+1}{p} \in (0,c+1)$, then the  families $(e^{zB^n_{m,p}})_{z \in \Sigma_\phi}$,  $(\sqrt z D_ye^{zB^n_{m,p}})_{z \in \Sigma_\phi}$, $\phi<\pi/2$ are $R$-bounded in $L^p_m(0,\infty)$. In particular,   $B^n_{m,p}$ has maximal regularity.
\end{teo}
{\sc Proof.} All the assertions follow from Theorem \ref{gen R-bound sempre m} and the heat kernel estimates of Propositions \ref{Estimates Bessel kernels}, \ref{Estimates gradient kernel bessel}, setting $\alpha=0$ or $\alpha=-1$ and $\beta=-c$. 
\qed

\begin{prop} \label{core2}
 If $c>-1$ and  $0<\frac{m+1}{p}< c+1$, then $$\mathcal D=\left \{v \in C_c^\infty (\R^{N}\times [0, \infty)), \ D_y v(x,0)=0 \right \}$$ is a core for $D(\mathcal L^n_{m,p})$.
\end{prop}
{\sc Proof.} Identical to that of Proposition \ref{core1}
\qed

\begin{teo}   \label{Bessel-estimates} Let  $c>-1$, $0<\frac{m+1}{p}<c+1$. Then $\mathcal L^n_{m,p}=\Delta_x+B^n_y$ with domain
$$D(\mathcal L^n_{m,p})=\{u \in W_m^{2,p}(\R^{N+1}_+):  \frac{D_y u}{y} \in L^p_m(\R^{N+1}_+)\} $$
has maximal regularity in $L^p_m(\R^{N+1}_+)$.
In particular, if $\frac{m+1}{p}<1$ then 
$$D(\mathcal L^n_{m,p})=\{u \in W_m^{2,p}(\R^{N+1}_+): D_y u(x,0)=0, \quad x \in \R^N \} $$ and 
if $\frac{m+1}{p}>1$ then 
$$D(\mathcal L^n_{m,p})= W_m^{2,p}(\R^{N+1}_+). $$

\end{teo}
{\sc Proof. } We apply Theorem \ref{Rbounded1} as in  Proposition \ref{analyt} and Theorem \ref{Domain} to prove that $\mathcal L^n_{m,p}$ has maximal regularity and is closed on the intersection of the domains of $\Delta_x$ and $B^n_y$. Then the same  argument  as in Theorem 
\ref{Lp-estimates} yield the $L^p_m$ boundedness of the mixed derivatives. By the domain characterization of $B^n$ and the closedness of $\Delta_x+B^n_y$ again we finally have  $D_{yy} u, y^{-1} D_y u \in  L^p_m (\R^{N+1}_+)$ for $u \in D(B^n_{m,p})$.

When $(m+1)/p>1$, by Hardy inequality of Proposition \ref{Hardy1} and  the equality $W^{1,p}_m=W^{1,p}_{0,m}$ of Proposition \ref{notrace},  we get $D(\mathcal L_{m,p})= W_m^{2,p}(\R^{N+1}_+) $. Instead, if $(m+1)/p<1$, then $(D_y u)/y \in L^p_m$ and $D_y u(x,0)=0$ are equivalent for $u \in W^{2,p}_m$. In fact, if $D_y u(x,0)=0$, then $D_y u \in W^{1,p}_{0,m}$ and $(D_y u)/y \in L^p_m$ by Proposition \ref{Hardy1}. Conversely, if $(D_y u)/y \in L^p_m$, since $\int_0^\infty \frac{|D_y u(x,y)|^p}{y^p}\, dy =\infty$ whenever $D_y u(x,0) \neq 0$, we get $D_y u(x,0)=0$ a.e.
\qed
Note that the closedness of $\Delta_x+B_y$ and the domain characterization of $B^n$ allow to conclude that $\Delta_x u, D_{yy} u \in L^p_m(\R^{N+1})$, for $u \in D(\mathcal L^n_{m,p})$, hence $D_{x_i x_j} u \in L^p_m(\R^{N+1})$ by the Calder\'on-Zygmund inequality in $\R^N$. However, to deduce that  the mixed derivatives $D_{x_i y} u$ belongs to $L^p_m(\R^{N+1})$ one needs the boundedness of singular integrals in $L^p_m(\R^{N+1})$ which usually requires that the one dimensional  weight $|y|^m$ belongs to $A_p(\R^{N+1})$, that is $0<(m+1)/p<1$, a restriction that does not appear in the above theorem.


\section{Rellich inequalities}
Rellich inequalities have been intensively studied in any dimension, also for degenerate operators, but here we recall only the 1d  result.
All the $L^p$ norms in this section are taken with respect to the Lebesgue measure and accordingly we write $L_p, \mathcal L_p$ for $L_{m,p}, \mathcal L_{m,p}$, when $m=0$.
We set for $1\le p \le \infty$
\begin{align*}
\gamma_p:&=\left (\frac{1}{p}-2\right )\left (\frac{1}{p'}+c\right )
\end{align*} 
and 
\begin{align} \label{Pp}
{\cal P}_{p}:&=\left\{\lambda=-\xi^2+i\xi \left (3-\frac{2}{p}+c \right )-\gamma_p\;;\;\xi\in \R\right\}.
\end{align}
Observe that ${\cal P}_{p}$ is a parabola with vertex at $-\gamma_p$ when $3-\frac{2}{p} +c \neq 0$, otherwise coincides with the semiaxis $]-\infty, -\gamma_p]$.

\begin{teo} \label{Rellich-p} (\cite[Section 3]{met-soba-spi},\cite[Section 4]{MNSS})
There exists a positive constant $C$ such that
\begin{equation} \label{rellich-op-y}
\left\||y|^{-2} u\right\|_p\leq C\| L u\|_p 
\end{equation}
holds for every $u\in D(L_{p,max})$ such that $u/|y|^2 \in L^p(0,\infty)$, if and only if,   $b \not \in {\cal P}_{p}$.
\end{teo}

Observe that $$b+\gamma_p=\left (\frac 1p -s_1-2 \right ) \left (s_2+2-\frac 1p \right )$$ so that $b \not \in {\cal P}_{p}$ means $\frac 1p \neq s_1+2, s_2+2$ when $3-\frac 2p+c \neq 0$ and $\frac 1p \in (s_1+2, s_2+2)$ when  $3-\frac 2p+c = 0$.

When  $s_1+2<\frac 1p < s_2+2$, independently of the value of $3-\frac 2p+c$, Rellich inequalities can be proved  by integrating by parts, see the Remark below. In the other cases the  proof relies on spectral theory and best constants are known only in special cases. We refer the reader to \cite[Section 3]{met-soba-spi}, again.

In the next result we show that Rellich inequalities hold for $\mathcal L=\Delta_x+L_y$ ($L_y=L$) in the generation interval $s_1<\frac 1p < s_2+2$  if and only if they hold for $L_y$. In the proof, the closedness of the sum $\Delta_x+L_y$ on the intersection of the corresponding domains plays a major role.

 Note that the theorem below (as that above) does not concern with Rellich inequalities in the whole domain of $\mathcal L$, as  described in the previous section, but on the (possibly) smaller subspace of all $u$ in the maximal domain, satisfying $u/|y|^2 \in L^p$.

\begin{teo} \label{Rellich-p-comp}
Assume $s_1<\frac 1p < s_2+2$. There exists $C>0$ such that
\begin{equation} \label{rellich-op}
\left\||y|^{-2} |u|\right\|_p\leq C\| \mathcal Lu\|_p 
\end{equation}
holds for every $u\in D(\mathcal L_{p,max})$  such that $u/|y|^2 \in L^p(\R_+^{N+1})$,  if and only if  $b \not \in {\cal P}_{p}$.
\end{teo}
{\sc Proof.} Assume that Rellich inequalities hold for the complete operator $\mathcal L$ and let $u(x,y)=z(x)\psi(y)$ with  $\|z\|_{L^p(\R^N)}=1$. Then 
$$\mathcal Lu=\Delta_x u+L_y u=\psi\Delta_x z+zL_y\psi$$ and
(\ref{rellich-op})  is equivalent to
\begin{equation*} 
\int_0^\infty\frac{|\psi(y)|^{p}}{|y|^{2p}}\, dy=\int_0^\infty|z(x)|^{p}\, dx\int_0^\infty\frac{|\psi(y)|^{p}}{|y|^{2p}}\, dy\leq C\int_{\R^{N+1}}\left (|\psi(y)\Delta_x z(x)+z(x) L_y\psi(y)\right |)^p\, dx\, dy.
\end{equation*}
Let $z_R(x)=R^{-\frac{N}{p}}z\left(\frac{x}{R}\right)$. Then $\|z_R\|_{L^p(\R^N)}=1$, $\|\Delta_x z_R\|_p\to 0$ as $R\to\infty$ and  letting $R \to \infty$
\begin{equation*} 
\int_0^\infty\frac{|\psi(y)|^{p}}{|y|^{2p}}\, dy\leq C\int_0^\infty|L_y\psi(y)|^p\, dy
\end{equation*}
so that Rellich inequalities hold for $L_y$.

Next, assume that Rellich inequalities hold for $L_y$ and let $u \in D(\mathcal L_{p,max})$ be such that $u/|y|^2 \in L^p(\R^{N+1})$. Then $u \in   D(\mathcal L_p)=D(\Delta_x)\cap D(L_{y,p})$ and for almost all $x \in \R^N$, $u(x, \cdot) \in D(L_{p,max})$ and $u(x, \cdot)/|y|^2 \in L^p((0, \infty))$. Then
$$
\int_0^\infty\frac{|u(x,y)|^p}{|y|^{2p}}\, dy \le  C\int_0^\infty|L_y u(x,y)|^p\, dy
$$
 and, integrating with respect to $x \in \R^N$ and using the closedness of $\Delta_x+L_y$  we get 
$$
\int_{\R^{N+1}}\frac{|u(x,y)|^p}{|y|^{2p}}\,dx\, dy \le  C\int_{\R^{N+1}}|L_y u(x,y)|^p\,dx\, dy \le C\int_{\R^{N+1}}|\mathcal Lu(x,y)|^p\, dx\, dy.
$$
 \qed

\begin{os} \label{Rellichparti}
When  $s_1+2<\frac 1p < s_2+2$ the best constant $C$ above is $(b+\gamma_p)^{-1}>0$. This can be seen   multiplying  $\mathcal L u$ by $u|u|^{p-2}y^{2-2p}$ and integrating by parts, assuming $u$ smooth and with support faraway from $\{y=0\}$). 
\end{os}

\section{Appendix A: auxiliary inequalities}
\begin{lem} \label{equiv}
For every $\eps>0$ there exists $C>0$ such that for $r,s >0$
\begin{equation*}
(1\wedge r)(1\wedge s) \le 1\wedge rs \le   C(1\wedge r)(1\wedge s)\,e^{\epsilon |r-s|^2}.
\end{equation*}
\end{lem} {\sc Proof.}  $(1\wedge rs) =(1\wedge r)(1\wedge s)$ when $r,s \le 1$ or $r,s \ge 1$. Assume that $s \le 1\le r$. Then $(1\wedge r)(1\wedge s)=s \le 1\wedge rs$. Conversely, if $rs \le 1$ then $1\wedge rs=rs \le Cse^{\eps(r-s)^2}$ for a suitable $C>0$, since $s \le 1\le r$. If, instead, $rs \ge 1$, then $1\wedge rs=1 \le Cr^{-1}e^{\eps(r-s)^2} \le Cse^{\eps(r-s)^2}$.
\qed

\begin{lem} \label{equiv1} If $\gamma_1 \le \gamma_2$ then for every $\eps>0$ there exists $C>0$ such that 
\begin{equation} \label{confrontoAlto}
\frac{|y|^{\gamma_1}}{|z|^{\gamma_2}}\leq C  \frac{(|y|\wedge 1)^{\gamma_1}}{(|z|\wedge 1)^{\gamma_2}} \exp\left(\epsilon |y-z|^2\right).
\end{equation}
\end{lem}
{\sc Proof. }
If $|y|\leq 1$ and $|z|\leq 1$ this is clearly true. 
Assume that $|z|\leq 1\leq |y|$. 
Then $|y-z|^2\geq (|y|-1)^2$ and
$$|y|^{\gamma_1}\leq C\exp\left\{\epsilon(|y|-1)^2\right\}\leq C\exp\left\{\epsilon(|y-z|)^2\right\}$$
and (\ref{confrontoAlto}) holds.
If $|y|\leq 1\leq |z|$ we  argue in similar way. Finally, when $|y| \ge 1$, $|z|\ge 1$ we write $y=r\omega, z=\rho \eta$ with $|\omega|=|\eta|=1$.
The left hand side of \eqref{confrontoAlto} is then $(r/\rho)^{\gamma_1}\rho^{\gamma_2-\gamma_1} \le (r/\rho)^{\gamma_1}$ which is now symmetric in $r, \rho$. 
Assuming  that $r \ge \rho \ge 1$ we write  $r=s \rho$ with $s \ge 1$, the inequality $s^{\gamma_1} \le C e^{\epsilon (s-1)^2 }   \le C e^{\epsilon (s-1)^2 \rho^2}$ ($s, \rho \ge 1$)  implies that 
$$
\frac{|y|^{\gamma_1}}{|z|^{\gamma_2}}\leq  \left (\frac{r}{\rho} \right )^{\gamma_1} \le C e^{\epsilon |r-\rho|^2} \le Ce^{ \epsilon |y-z|^2}.
$$
\qed

The following Hardy inequalities have been used several times throughout the paper.
\begin{lem} \label{Hardy} 
\begin{itemize}
\item[(i)] When $c+1>\frac{m+1}{p}$ the map $H_1f(y)=\frac{1}{y^{c+1}} \int_0^y f(s) s^c\, ds$ is bounded from $L^p_m$ to itself. \\
\item[(ii)] When $c+1<\frac{m+1}{p}$ the map $H_2f(y)=\frac{1}{y^{c+1}} \int_y^\infty f(s) s^c\, ds$ is bounded from $L^p_m$ to itself.
\end{itemize}
\end{lem}
{\sc Proof.} Take $f \ge 0$ and let $w=H_1f$. Then $w(y)=\int_0^1 f(ty)t^c\, dt$ and by Minkowski's inequality
$$
\|w\|_{L^p_m} \le \int_0^1 t^c\left ( \int_0^\infty f(ty)^py^mdy \right )^{\frac1p}=\int_0^1 t^{c-\frac{m+1}{p}}\left ( \int_0^\infty f(x)^py^mdx \right )^{\frac1p}=C\|f\|_{L^p_m}
$$
with $C=\left((c+1)-\frac{m+1}{p}\right)^{-1}$. The proof for $H_2$ is similar.
\qed

\section{Appendix B: Sobolev spaces with weights}\label{Appendix Sobolev}
For $1<p<\infty$ let $W^{k,p}_m(\R^{N+1}_+)=\{u \in L^p_m(\R^{N+1}_+): \partial^\alpha u \in  L^p_m(\R^{N+1}_+) \quad |\alpha| \le k\}$. We use often $W^{k,p}_m$ thus omitting $\R^{N+1}_+$ and $W^{k,p}_{0,m}$ for the closure of $C_c^\infty (\R^{N+1}_+)$ in $W^{k,p}_m$.
 \begin{lem} \label{l1}
If $\frac{m+1}{p}<1$ then $L^p_m$ embeds into $L^1(Q \times (0,1))$, where $Q$ is any cube of $\R^N$. It follows that $W^{1,p}_m$ embeds into $W^{1,1}(Q\times (0,1))$ and that every function $u\in W^{1,p}_m$ has a trace $u(\cdot, y) \in L^1(Q)$, for every $0 \le y \le 1$.
\end{lem}
{\sc Proof. }
$$\int_{Q\times (0,1)} |u|dxdy=\int_{Q\times (0,1)} |u|y^{\frac mp} y^{- \frac mp} dxdy \le \|u\|_{L^p_m} \left (\int_{Q \times (0,1 )} y^{-\frac{m}{p-1}} dx dy\right)^{\frac{1}{p'}}.
$$
\qed
\smallskip
Keeping the assumption $(m+1)/p<1$ we write for $u \in W^{1,p}_m$ (use that for almost all $x \in \R^N$, $u(x, \cdot) \in W^{1,1}(0,1)$)
$$
|u(x,y)-u(x,0)|=\left|\int_0^y D_y u(x,s)s^{\frac mp} s^{-\frac mp}\, ds\right| \le C\left (\int_0^y |D_y u(x,s)|^p s^m ds\right )^{\frac 1p} y^{1-\frac{m+1}{p}}.
$$
Raising to the power $p$ and integrating with respect to $x$ we get
\begin{equation} \label{fond}
\int_{\R^N} |u(x,y)-u(x,0)|^p  dx\le  Cy^{p-m-1} \int_{\R^N \times (0,y) }|D_y u(x,s)|^p s^m dx ds, \quad \frac{m+1}{p} <1.
\end{equation}
\begin{lem} \label{mle-1}
If $m \le -1$ and $u \in W^{1,p}_m$, then $u(\cdot,0)=0$.
\end{lem}
{\sc Proof. }If $E=\{x \in \R^N: u(x,0) \neq 0 \}$, then 
$
\int_0^1 |u(x,y)|^py^m  dy =\infty$
for every $x \in E$, since $m \le -1$. But then $|E|=0$.
\qed

\begin{prop} \label{trace} If $\frac{m+1}{p} <1$, then $u \in W^{1,p}_{0,m}$ if and only if $u(\cdot, 0)=0$. In particular $W^{1,p}_{m}=W^{1,p}_{0,m}$ if $m \le -1$.
\end{prop}
{\sc Proof. } If $u \in W^{1,p}_{0,m}$, then $u \in W^{1,1}_0\left(Q\times (0,1)\right)$ and then $u(\cdot, 0)=0$. To show the converse, it suffices to show that $u$ can be approximated in the $W^{1,p}_m$ norm by functions with support far away from $y=0$. After this, cut-off near infinity and smoothing using convolutions is standard. Let $\phi$ be a smooth function which is equal to $0$ in $(0,1)$ and to $1$ for $y \ge 2$ and $u_n(x,y)=\phi(ny)u(x,y)$. By dominated convergence $u_n \to u, D_{x_i} u_n \to D_{x_i} u$ in $L^p_m$. Moreover, $D_y u_n =\phi (ny) D_y u+n \phi'(ny) u$ and $\phi(ny)D_y u \to D_y u$, by dominated convergence, again. Finally, multiplying  \eqref{fond} by $y^m$ and integrating with respect to $y$ we get
\begin{eqnarray*}
\int_{\R^{N+1}_+}n^p |\phi'(ny)|^p |u(x,y)|^p y^m dxdy &\le& \|\phi'\|_\infty^p n^p \int_{\R^N \times (0, \frac 2n)}|u(x,y)|^p y^m dx dy  \\[1ex]
&\le & C\|\phi'\|_\infty^p  \int_{\R^N \times (0, \frac 2n)}|D_y u(x,y)|^p y^m dx dy \to 0 \quad {\rm as}\ n \to \infty.
\end{eqnarray*}
\qed
\begin{prop} \label{notrace} If $\frac{m+1}{p} \ge1$, then  $W^{1,p}_{m}=W^{1,p}_{0,m}$.
\end{prop}
{\sc Proof. } As in the proof above, we approximate with functions whose support is far away $y=0$, the rest being standard. We may also suppose that $u$ is bounded (otherwise we consider $(u\wedge n) \vee (-n))$ and with  support contained in $Q \times (0, \ell)$, $\ell \ge 1$.  Let $\phi$ be a smooth function which is equal to $0$ in $\left(0,\frac{1}{4}\right)$ and to $1$ for $y \ge \frac{1}{2}$. In order to consider also the critical case $\frac{m+1}{p}=1$, we choose $\phi_n(y)=\phi\left(y^\frac{1}{n}\right)$ and set $u_n(x,y)=\phi\left(y^\frac{1}{n}\right)u(x,y)$. By dominated convergence $u_n \to u, D_{x_i} u_n \to D_{x_i} u$ in $L^p_m$. Moreover, $D_y u_n =\phi_n (y) D_y u+\frac{1}{n}\phi'\left(y^\frac{1}{n}\right)y^{\frac{1}{n}-1}  u$ and $\phi_n(y)D_y u \to D_y u$, by dominated convergence, again. It remains  to show that $\frac{1}{n}\phi'\left(y^\frac{1}{n}\right)y^{\frac{1}{n}-1} u \to 0$ in $L^p_m$.  Using $m+1-p\geq 0$ then  we get 
\begin{align*}
&\frac{1}{n^p}\int_{Q\times (0,\ell)}  \left|\phi'\left(y^\frac{1}{n}\right)\right|^p |u(x,y)|^p y^{m+p(\frac1n-1)} dxdy\leq 
\frac 1{n^p}\|u\|^p_\infty \int_{Q\times \left((\frac 1 4)^n,\,(\frac 1 2)^n \right)}|\phi'(y^\frac 1n)|^py^{p\frac 1 n-p+m}\,dy \\[1ex]
&=\frac1{n^{p-1}}\|u\|^p_\infty |Q|\int_{\frac 1 4}^{\frac 1 2}|\phi'(s)|^ps^{p-1+n(m+1-p)}\,ds\leq \frac1{n^{p-1}}\|u\|^p_\infty |Q|\int_{\frac 1 4}^{\frac 1 2}|\phi'(s)|^p\,ds\to 0
\end{align*}
 as $n$ goes to infinity, since $p>1$.
\qed

\begin{os} If $(m+1)/p \geq 1$ and $p>1$ it is non true, in general, that $u\in W_m^{1,p}$ has a finite trace on the boundary. For example, if 
$m+1=p$ and  $N=0$  one can take 
 $u(y)=\log(\log y^{-1})$; if $m+1>p$ one can take $y^{1-a}$ with $1<a<\frac{m+1}p$.    
\end {os}
 

 Hardy inequalities in $W^{1,p}_{0,m}$ follow from Lemma \ref{Hardy}.
\begin{prop} \label{Hardy1} If $\frac{m+1}{p} \neq 1$, then Hardy inequality
$$
\left\|\frac{u}{y}\right\|_{L^p_m} \le C\|D_y u\|_{L^p_m}
$$ holds for $u \in W^{1,p}_{0,m}$.
\end{prop}
{\sc Proof. } This follows from Lemma \ref{Hardy} (i), (ii) with $c=0$, $f =D_yu $, after integrating with respect to $x \in \R^N$ the $1d$ inequalities.
\qed

\section {Appendix C: boundedness of a family of integral operators}

Next we consider a two-parameters  family of integral operators $\left(S_{\alpha,\beta}(t)\right)_{t>0}$ $L^p_m(\R^M)$, defined for $\alpha,\beta\in\R$  and $t>0$ by
\begin{align*}
S^{\alpha,\beta}(t)f(y)=t^{-\frac {M} 2}\,\left (\frac{|y|}{\sqrt t}\wedge 1 \right)^{-\alpha}\int_{\R^M}  \left (\frac{|z|}{\sqrt t}\wedge 1 \right)^{-\beta}
\exp\left(-\frac{|y-z|^2}{\kappa t}\right)f(z) \,dz,
\end{align*}
where $\kappa$ is a positive constant. We omit the dependence on $\kappa$ even though in some proofs we need to vary it.

We start by observing  that the  scale homogeneity of $S^{\alpha,\beta}$ is $2$ since a simple change of variable in the integral yields
\begin{align*}
S^{\alpha,\beta}(t)\left(I_s f\right)=I_s \left(S^{\alpha,\beta}(s^2t) f\right),\qquad I_sf(y)=f(sy),\qquad t, s>0,
\end{align*}
which in particular gives 
\begin{align*}
S^{\alpha,\beta}(t)f =I_{1/\sqrt t} \left(S^{\alpha,\beta}(1) I_{\sqrt t} f\right),\qquad t> 0.
\end{align*}
The boundedness of $S^{\alpha,\beta}(t)f $ in $L^p_m(\R^M)$ is then equivalent to that for $t=1$ and $\|S^{\alpha,\beta}(t)\|_p= \|S^{\alpha,\beta}(1)\|_p$.

Let us state first a simple lemma which allows to change the reference measure without additional efforts.
\begin{lem}\label{Lem Equiv integral oper}
For  $\alpha,\beta,\gamma_1 \le \gamma_2	\in\R$ let 
\begin{align}\label{Equiv integral oper}
\left(\tilde S^{\alpha,\beta,\gamma_1, \gamma_2}(t)\right)f(y)=t^{-\frac {M} 2}\,\left (\frac{|y|}{\sqrt t}\wedge 1 \right)^{-\alpha}\int_{\R^M}\frac{|y|^{\gamma_1}}{|z|^{\gamma_2}}  \left (\frac{|z|}{\sqrt t}\wedge 1 \right)^{-\beta}
\exp\left(-\frac{|y-z|^2}{\kappa t}\right)f(z) \,dz.
\end{align}
Then for every $\kappa'>\kappa$  and $f \ge 0$
\begin{align*}
\left(\tilde S^{\alpha,\beta,\gamma_1, \gamma_2}(t)\right)f(y)&\leq C t^{-\frac {M+\gamma_1-\gamma_2} 2}\,\left (\frac{|y|}{\sqrt t}\wedge 1 \right)^{-\alpha+\gamma_1}\int_{\R^M}\left (\frac{|z|}{\sqrt t}\wedge 1 \right)^{-\beta-\gamma_2}
\exp\left(-\frac{|y-z|^2}{\kappa' t}\right)f(z) \,dz\\[1ex]
&=Ct^{-\frac{\gamma_1-\gamma_2}{2}}S^{\alpha-\gamma_1,\beta+\gamma_2}(t)f(y).
\end{align*}
\end{lem}
{\sc{Proof.}} This follows from Lemma \ref{equiv1}.
\qed


\begin{prop}\label{Boundedness theta}
Let $m\in\R$, $\theta\geq 0$. The following properties are equivalent.
\begin{itemize}
\item[(i)]  For every $t>0$ $(S^{\alpha,\beta}(t))_{t \ge 0}$ is  bounded from   $L^p_m$ to $ L^p_{m-p\theta}$ and
\begin{align*}
\|S^{\alpha,\beta}(t)\|_{L^p_m\to L^p_{m-p\theta}}\leq Ct^{-\frac\theta 2}.
\end{align*}
\item[(ii)] $S^{\alpha,\beta}(1)$ is  bounded from  from   $L^p_m$ to $ L^p_{m-p\theta}$.
\item[(iii)] $\alpha<\frac{M+m}p < M-\beta$ (in  particular $\alpha+\beta<M-\theta$).
\end{itemize}
\end{prop}
{\sc Proof.} The equivalence of (i) and (ii) follows by using a scaling argument since
\begin{align*}
S^{\alpha,\beta}(t)f =I_{1/\sqrt t} \left(S^{\alpha,\beta}(1) I_{\sqrt t} f\right),\qquad t> 0,
\end{align*}
where $I_sf(y)=f(sy)$, $s>0$ and $\|I_s\|_{L^p_m}=s^{-\frac{M+m}p}$.
The boundedness of $S^{\alpha,\beta}(t)f $ from   $L^p_m$ to $ L^p_{m-p\theta}$ is then equivalent to that for $t=1$ and 
$$\|S^{\alpha,\beta}(t)\|_{L^p_m\to L^p_{m-p\theta}}= \|S^{\alpha,\beta}(1)\|_{L^p_m\to L^p_{m-p\theta}}=t^{-\frac \theta 2}$$

 Condition (iii) means that  $|y|^{-\alpha}\in L^p_{m-p\theta}(B, d\mu_m)$, $ |y|^{-\beta-m}\in L^{p'}_{m}(B)$, where $B$ is the unit ball, and these are necessary for the boundedness of $S^{\alpha, \beta}(1)$.
 
  To prove that (iii) implies (ii), let us assume first that $m=0$ and write $S^{\alpha, \beta}(1)f=S_1(f)+S_2(f)$ where for $f \ge 0$
\begin{align*}
S_1(f)(y)& =(|y|\wedge 1)^{-\alpha} \int_B |z|^{-\beta} e^{-\frac{|y-z|^2}{\kappa}} f(z)dz \le C(|y|\wedge 1)^{-\alpha}e^{-c|y|^2} \int_B |z|^{-\beta}f(z)dz \\
&\le C\|f\|_p 
(|y|\wedge 1)^{-\alpha} e^{-c|y|^2},
\end{align*}
(using that $|z| \le 1$) and 
\begin{align*}
S_2(f)(y)=(|y|\wedge 1)^{-\alpha} \int_{\R^M \setminus B}  e^{-\frac{|y-z|^2}{\kappa}}f(z)dz. 
\end{align*}
The boundedness of $S_1$ is clear. The one  of $S_2$  for $|y| \ge 1$  follows after observing that in this range $|y|^{-p\theta}\leq 1$ which   implies $\|S_2\|_{L^p_{-p\theta}(B^c)}\leq \|S_2\|_{L^p(B^c)}$: the result then follows by using Young's inequality. On the other hand the required estimates for $S_2$ in the range   $|y| \le 1$ follows by estimating the sup-norm of the convolution with the $L^p$ norm of $f$.

For a generic $m$ we consider the isometries 
\begin{align*}
T_{-\frac m p}&: L^p\to L^p_m\quad \hspace{5.7ex}f\mapsto |y|^{-\frac m p}f,\\[1ex]
T_{\frac m p}&: L^p_{m-p\theta}\to L^p_{-p\theta}\quad f\mapsto |y|^{\frac m p}f.
 \end{align*}
Then,  using Lemma \ref{Lem Equiv integral oper},
\begin{align*}
|T_{\frac m p} S^{\alpha,\beta}(1)T_{-\frac m p}f|=|\tilde S^{\alpha,\beta,\frac{m}{p}, \frac{m}{p}}(1)f|\leq  C S^{\alpha-\frac m p,\beta+\frac{m}p}(1)|f|
\end{align*}
and the thesis follows from the case $m=0$.
\qed

\end{document}